# $C^*$-Rhapsody

Aleks Kleyn

Abstract. In the paper I considered linear and antilinear mappings of finite dimensional algebra over complex field, as well I considered involution. I considered also an example of $C^*$-algebra.

## Contents

## 1. Preface

Since I have been studying calculus in Banach algebra, it is difficult to ignore $C^*$-algebra which is a tool in quantum mechanics. However, before considering particular statements, I need to understand the structure of the linear mapping of $C^*$-algebra. This caused to write this paper.

However there exist different algebras where conjugation is defined. For this reason, I decided not restrict myself by considering of complex field. Considering of linear and antilinear mappings of algebra with conjugation turned out to be different from the concept, which I considered in [4]. It is attractive to find a statement similar to statement that we can expand additive mapping of complex field into sum of linear and antilinear mappings. Such statement would be a bridge

Aleks_Kleyn@MailAPS.org.
http://sites.google.com/site/AleksKleyn/.
http://arxiv.org/a/kleyn_a_1.
http://AleksKleyn.blogspot.com/.

between the results in [4] and this paper. However, there is reason to believe that in general this is not true.

This means that, in general, consideration of the conjugation is not enough for the complete classification of linear maps. However, this classification is important to understand the structure of closed differential form in Banach algebra.

## 2. Conventions

(1) Let $A$ be free finite dimensional algebra. Considering expansion of element of algebra $A$ relative basis $\overline{\overline{e}}$ we use the same root letter to denote this element and its coordinates. However we do not use vector notation in algebra. In expression $a^2$, it is not clear whether this is component of expansion of element $a$ relative basis, or this is operation $a^2 = aa$. To make text clearer we use separate color for index of element of algebra. For instance,
$$a = a^{i}\overline{e}_{i}$$
(2) If free finite dimensional algebra has unit, then we identify the vector of basis $\overline{e}_{0}$ with unit of algebra.
(3) Although the algebra is a free module over some ring, we do not use the vector notation to write elements of algebra. In the case when I consider the matrix of coordinates of element of algebra, I will use vector notation to write corresponding element. In order to avoid ambiguity when I use conjugation, I denote $a^*$ element conjugated to element $a$.
(4) Without a doubt, the reader may have questions, comments, objections. I will appreciate any response.

## 3. $C$-algebra

Considered below construction is based on the construction considered in section [3]-5. In this paper, we will use notation offered in this section.

Let $C$ be complex field. Let $\overline{\overline{e}}_C$

$$\overline{e}_{C\cdot 0} = 1 \quad \overline{e}_{C\cdot 1} = i \tag{3.1}$$

be the basis of algebra $C$ over real field $R$. Product in algebra $C$ is defined according to rule

$$\overline{e}^2_{C\cdot 1} = -\overline{e}_{C\cdot 0} \tag{3.2}$$

According to the theorem [3]-3.1, structural constants of complex field $C$ over real field $R$ have form

$$\begin{array}{ll} C_{C\cdot}{}^{0}_{00} = 1 & C_{C\cdot}{}^{1}_{01} = 1 \\ C_{C\cdot}{}^{1}_{10} = 1 & C_{C\cdot}{}^{0}_{11} = -1 \end{array} \tag{3.3}$$

According to definition [4]-2.2.1, $C$-algebra $A$ is $C$-vector space. Let $\overline{\overline{e}}_{AC}$ be a $C$-basis of finite dimensional $C$-algebra $A$. The product in $C$-algebra $A$ is defined by equation

$$\overline{e}_{AC\cdot k}\overline{e}_{AC\cdot l} = C_{AC\cdot}{}^{i}_{kl}\overline{e}_{AC\cdot i} \tag{3.4}$$

where $C_{AC\cdot}{}^{k}_{ij}$ are structural constants of $C$-algebra $A$ over the field $C$.

I will consider the $C$-algebra $A$ as direct sum of algebras $C$. Each item of sum I identify with vector of basis $\overline{\overline{e}}_{AC}$ . Accordingly, I can consider $C$-algebra $A$ as algebra over field $R$. Let $\overline{\overline{e}}_{AR}$ be basis of $C$-algebra $A$ over the field $R$. Index of basis $\overline{\overline{e}}_{AR}$ consists from two indexes: index of fiber and index of vector of basis $\overline{\overline{e}}_C$ in fiber.

I will identify vector of basis $\overline{e}_{AC\cdot i}$ with unit in corresponding fiber. Then

$$\overline{e}_{AR\cdot ji}=\overline{e}_{C\cdot j}\overline{e}_{AC\cdot i} \tag{3.5}$$

The product of vectors of basis $\overline{\overline{e}}_{AR}$ has form

$$\overline{e}_{AR\cdot ji}\overline{e}_{AR\cdot mk}=\overline{e}_{C\cdot j}\overline{e}_{AC\cdot i}\overline{e}_{C\cdot m}\overline{e}_{AC\cdot k}=C_{C\cdot}{}^{a}_{jm}\overline{e}_{C\cdot a}C_{AC\cdot}{}^{b}_{ik}\overline{e}_{AC\cdot b} \tag{3.6}$$

$$\overline{e}_{AR\cdot 00}\,\overline{e}_{AR\cdot ki}=\overline{e}_{AR\cdot ki} \tag{3.7}$$

$$\overline{e}_{AR\cdot 10}\,\overline{e}_{AR\cdot 1i}=-\overline{e}_{AR\cdot 0i} \tag{3.8}$$

Because $C_{AC\cdot}{}^{b}_{ik}\in C$, then expansion $C_{AC\cdot}{}^{b}_{ik}\in C$ relative to basis $\overline{\overline{e}}_C$ has form

$$C_{AC\cdot}{}^{b}_{ik}=C_{AC\cdot}{}^{bc}_{ik}\overline{e}_{C\cdot c} \tag{3.9}$$

We can represent transformation $f$ using matrix $f^{ij}_{kl}$

$$a'^{ij}=f^{ij}_{kl}a^{kl} \tag{3.10}$$

$$\begin{aligned}a'^{0j}&=f^{0j}_{0l}a^{0l}+f^{0j}_{1l}a^{1l}\\ a'^{1j}&=f^{1j}_{0l}a^{0l}+f^{1j}_{1l}a^{1l}\end{aligned} \tag{3.11}$$

We can represent $c\in C$ in following form

$$c=c^{0}\,\overline{e}_{AR\cdot 00}+c^{1}\,\overline{e}_{AR\cdot 10} \tag{3.12}$$

$$c^{*}=c^{0}\,\overline{e}_{AR\cdot 00}-c^{1}\,\overline{e}_{AR\cdot 10} \tag{3.13}$$

From equations (3.7), (3.8), (3.12), it follows

$$(ca)^{0i}=c^{0}a^{0i}-c^{1}a^{1i} \tag{3.14}$$

$$(ca)^{1i}=c^{0}a^{1i}+c^{1}a^{0i} \tag{3.15}$$

**Theorem 3.1.** *Structural constants of $C$-algebra $A$ over real field $R$ have form*

$$\left\{\begin{aligned}
C_{AR\cdot}{}^{\cdot 0b}_{\cdot 0i\cdot 0k}&=C_{AC\cdot}{}^{b0}_{ik}\\
C_{AR\cdot}{}^{\cdot 0b}_{\cdot 1i\cdot 1k}&=-C_{AC\cdot}{}^{b0}_{ik}\\
C_{AR\cdot}{}^{\cdot 0b}_{\cdot 0i\cdot 1k}&=-C_{AC\cdot}{}^{b1}_{ik}\\
C_{AR\cdot}{}^{\cdot 0b}_{\cdot 1i\cdot 0k}&=-C_{AC\cdot}{}^{b1}_{ik}\\
C_{AR\cdot}{}^{\cdot 1b}_{\cdot 0i\cdot 1k}&=C_{AC\cdot}{}^{b0}_{ik}\\
C_{AR\cdot}{}^{\cdot 1b}_{\cdot 1i\cdot 0k}&=C_{AC\cdot}{}^{b0}_{ik}\\
C_{AR\cdot}{}^{\cdot 1b}_{\cdot 0i\cdot 0k}&=C_{AC\cdot}{}^{b1}_{ik}\\
C_{AR\cdot}{}^{\cdot 1b}_{\cdot 1i\cdot 1k}&=-C_{AC\cdot}{}^{b1}_{ik}
\end{aligned}\right. \tag{3.16}$$

*Proof.* According to equations [3]-(5.3), [3]-(5.5), structural constants of algebra $A$ over field $R$ have form

$$\begin{aligned}C_{AR\cdot}{}^{\cdot db}_{\cdot ji\cdot mk}&=C_{C\cdot}{}^{a}_{jm}C_{C\cdot}{}^{d}_{ac}C_{AC\cdot}{}^{bc}_{ik}\\ C_{AC\cdot}{}^{b}_{ik}&=C_{AC\cdot}{}^{bc}_{ik}\overline{e}_{C\cdot c}=C_{AC\cdot}{}^{b0}_{ik}+C_{AC\cdot}{}^{b1}_{ik}i\end{aligned} \tag{3.17}$$

From equations (3.3), (3.17), it follows that

$$
\begin{cases}
C_{AR\cdot}{}^{\cdot 0b}_{ji\cdot mk} = C_{C\cdot}{}^{0}_{jm} C_{C\cdot}{}^{0}_{00} C_{AC\cdot}{}^{b0}_{ik} + C_{C\cdot}{}^{1}_{jm} C_{C\cdot}{}^{0}_{11} C_{AC\cdot}{}^{b1}_{ik}\\
\qquad = C_{C\cdot}{}^{0}_{jm} C_{AC\cdot}{}^{b0}_{ik} - C_{C\cdot}{}^{1}_{jm} C_{AC\cdot}{}^{b1}_{ik}\\
C_{AR\cdot}{}^{\cdot 1b}_{ji\cdot mk} = C_{C\cdot}{}^{1}_{jm} C_{C\cdot}{}^{1}_{10} C_{AC\cdot}{}^{b0}_{ik} + C_{C\cdot}{}^{0}_{jm} C_{C\cdot}{}^{1}_{01} C_{AC\cdot}{}^{b1}_{ik}\\
\qquad = C_{C\cdot}{}^{1}_{jm} C_{AC\cdot}{}^{b0}_{ik} + C_{C\cdot}{}^{0}_{jm} C_{AC\cdot}{}^{b1}_{ik}
\end{cases}
\tag{3.18}
$$

From equations (3.3), (3.18), it follows that

$$
\begin{cases}
C_{AR\cdot}{}^{\cdot 0b}_{0i\cdot 0k} = C_{C\cdot}{}^{0}_{00} C_{AC\cdot}{}^{b0}_{ik} - C_{C\cdot}{}^{1}_{00} C_{AC\cdot}{}^{b1}_{ik}\\
C_{AR\cdot}{}^{\cdot 0b}_{1i\cdot 1k} = C_{C\cdot}{}^{0}_{11} C_{AC\cdot}{}^{b0}_{ik} - C_{C\cdot}{}^{1}_{11} C_{AC\cdot}{}^{b1}_{ik}\\
C_{AR\cdot}{}^{\cdot 0b}_{0i\cdot 1k} = C_{C\cdot}{}^{0}_{01} C_{AC\cdot}{}^{b0}_{ik} - C_{C\cdot}{}^{1}_{01} C_{AC\cdot}{}^{b1}_{ik}\\
C_{AR\cdot}{}^{\cdot 0b}_{1i\cdot 0k} = C_{C\cdot}{}^{0}_{10} C_{AC\cdot}{}^{b0}_{ik} - C_{C\cdot}{}^{1}_{10} C_{AC\cdot}{}^{b1}_{ik}\\
C_{AR\cdot}{}^{\cdot 1b}_{0i\cdot 1k} = C_{C\cdot}{}^{1}_{01} C_{AC\cdot}{}^{b0}_{ik} + C_{C\cdot}{}^{0}_{01} C_{AC\cdot}{}^{b1}_{ik}\\
C_{AR\cdot}{}^{\cdot 1b}_{1i\cdot 0k} = C_{C\cdot}{}^{1}_{10} C_{AC\cdot}{}^{b0}_{ik} + C_{C\cdot}{}^{0}_{10} C_{AC\cdot}{}^{b1}_{ik}\\
C_{AR\cdot}{}^{\cdot 1b}_{0i\cdot 0k} = C_{C\cdot}{}^{1}_{00} C_{AC\cdot}{}^{b0}_{ik} + C_{C\cdot}{}^{0}_{00} C_{AC\cdot}{}^{b1}_{ik}\\
C_{AR\cdot}{}^{\cdot 1b}_{1i\cdot 1k} = C_{C\cdot}{}^{1}_{11} C_{AC\cdot}{}^{b0}_{ik} + C_{C\cdot}{}^{0}_{11} C_{AC\cdot}{}^{b1}_{ik}
\end{cases}
\tag{3.19}
$$

(3.16) follows from equations (3.3), (3.19). □

From the equation (3.16), it follows that

$$
\begin{cases}
C_{AR\cdot}{}^{\cdot 0b}_{0i\cdot 0k} = -C_{AR\cdot}{}^{\cdot 0b}_{1i\cdot 1k} = C_{AC\cdot}{}^{b0}_{ik}\\
C_{AR\cdot}{}^{\cdot 0b}_{0i\cdot 1k} = C_{AR\cdot}{}^{\cdot 0b}_{1i\cdot 0k} = -C_{AC\cdot}{}^{b1}_{ik}\\
C_{AR\cdot}{}^{\cdot 1b}_{0i\cdot 1k} = C_{AR\cdot}{}^{\cdot 1b}_{1i\cdot 0k} = C_{AC\cdot}{}^{b0}_{ik}\\
C_{AR\cdot}{}^{\cdot 1b}_{0i\cdot 0k} = -C_{AR\cdot}{}^{\cdot 1b}_{1i\cdot 1k} = C_{AC\cdot}{}^{b1}_{ik}
\end{cases}
\tag{3.20}
$$

## 4. Linear Mapping of $C$-algebra

**Definition 4.1.** The mapping of $C$-vector spaces

$$f : A_1 \to A_2$$

is called $C$-linear if

$$f(ca) = cf(a) \quad c \in C \quad a \in A \tag{4.1}$$

□

**Theorem 4.2** (The Cauchy-Riemann equations)**.** *Matrix of linear mapping of $C$-vector spaces*

$$f : A_1 \to A_2$$
$$y^{ik} = f^{ik}_{jm} x^{jm}$$

*satisfies relationship*

$$
\begin{aligned}
f^{\cdot 0j}_{\cdot 0i} &= f^{\cdot 1j}_{\cdot 1i}\\
f^{\cdot 1j}_{\cdot 0i} &= -f^{\cdot 0j}_{\cdot 1i}
\end{aligned}
\tag{4.2}
$$

*Proof.* From equations (3.11), (3.14), (3.15), it follows that

$$\begin{aligned}(f(ca))^{0j} &= f^{0j}_{0i}(ca)^{0i} + f^{0j}_{1i}(ca)^{1i}\\ &= f^{0j}_{0i}(c^0a^{0i} - c^1a^{1i}) + f^{0j}_{1i}(c^0a^{1i} + c^1a^{0i})\\ &= f^{0j}_{0i}c^0a^{0i} - f^{0j}_{0i}c^1a^{1i} + f^{0j}_{1i}c^0a^{1i} + f^{0j}_{1i}c^1a^{0i}\end{aligned} \tag{4.3}$$

$$\begin{aligned}(f(ca))^{1j} &= f^{1j}_{0i}(ca)^{0i} + f^{1j}_{1i}(ca)^{1i}\\ &= f^{1j}_{0i}(c^0a^{0i} - c^1a^{1i}) + f^{1j}_{1i}(c^0a^{1i} + c^1a^{0i})\\ &= f^{1j}_{0i}c^0a^{0i} - f^{1j}_{0i}c^1a^{1i} + f^{1j}_{1i}c^0a^{1i} + f^{1j}_{1i}c^1a^{0i}\end{aligned} \tag{4.4}$$

From equations (4.1), (3.11), (3.12), (3.14), (3.15), it follows that

$$\begin{aligned}(cf(a))^{0j} &= c^0(f(a))^{0j} - c^1(f(a))^{1j}\\ &= c^0(f^{0j}_{0i}a^{0i} + f^{0j}_{1i}a^{1i}) - c^1(f^{1j}_{0i}a^{0i} + f^{1j}_{1i}a^{1i})\\ &= c^0f^{0j}_{0i}a^{0i} + c^0f^{0j}_{1i}a^{1i} - c^1f^{1j}_{0i}a^{0i} - c^1f^{1j}_{1i}a^{1i}\end{aligned} \tag{4.5}$$

$$\begin{aligned}(cf(a))^{1j} &= c^0(f(a))^{1j} + c^1(f(a))^{0j}\\ &= c^0(f^{1j}_{0i}a^{0i} + f^{1j}_{1i}a^{1i}) + c^1(f^{0j}_{0i}a^{0i} + f^{0j}_{1i}a^{1i})\\ &= c^0f^{1j}_{0i}a^{0i} + c^0f^{1j}_{1i}a^{1i} + c^1f^{0j}_{0i}a^{0i} + c^1f^{0j}_{1i}a^{1i}\end{aligned} \tag{4.6}$$

From equations (4.1), (4.3), (4.5), it follows that

$$\begin{aligned}&f^{0j}_{0i}c^0a^{0i} - f^{0j}_{0i}c^1a^{1i} + f^{0j}_{1i}c^0a^{1i} + f^{0j}_{1i}c^1a^{0i}\\ =&c^0f^{0j}_{0i}a^{0i} + c^0f^{0j}_{1i}a^{1i} - c^1f^{1j}_{0i}a^{0i} - c^1f^{1j}_{1i}a^{1i}\end{aligned} \tag{4.7}$$

From equations (4.1), (4.4), (4.6), it follows that

$$\begin{aligned}&f^{1j}_{0i}c^0a^{0i} - f^{1j}_{0i}c^1a^{1i} + f^{1j}_{1i}c^0a^{1i} + f^{1j}_{1i}c^1a^{0i}\\ =&c^0f^{1j}_{0i}a^{0i} + c^0f^{1j}_{1i}a^{1i} + c^1f^{0j}_{0i}a^{0i} + c^1f^{0j}_{1i}a^{1i}\end{aligned} \tag{4.8}$$

Since $c^k$, $a^{ki}$ are arbitrary, then the equation (4.2) follows from equations (4.7), (4.8). □

**Theorem 4.3.** *Let $\overline{\overline{e}}_{iC}$ be C-basis of C-vector space $A_i$, $i = 1, 2$. Let $\overline{\overline{e}}_{iR}$ be corresponding R-basis of C-vector space $A_i$. The matrix of linear mapping*

$$f : A_1 \to A_2$$

*has form*

$$\begin{aligned}x^{0j}_2 &= f^{0j}_{0i}x^{0i}_1 - f^{1j}_{0i}x^{1i}_1\\ x^{1j}_2 &= f^{1j}_{0i}x^{0i}_1 + f^{0j}_{0i}x^{1i}_1\end{aligned} \tag{4.9}$$

*relative bases $\overline{\overline{e}}_{1R}$, $\overline{\overline{e}}_{2R}$ and has form*

$$x^j_2 = f^j_i x^i_1 \quad f^j_i = f^{0j}_{0i} + f^{1j}_{0i}i \tag{4.10}$$

*relative bases $\overline{\overline{e}}_{1C}$, $\overline{\overline{e}}_{2C}$*

*Proof.* Since the mapping $f$ is $R$-linear mapping, then relative basis $\overline{\overline{e}}_{1R}$, $\overline{\overline{e}}_{2R}$ mapping $f$ has form

$$x_2^{\mathbf{0}j} = f_{\mathbf{0}i}^{\mathbf{0}j} x_1^{\mathbf{0}i} + f_{\mathbf{1}i}^{\mathbf{0}j} x_1^{\mathbf{1}i}$$
$$x_2^{\mathbf{1}j} = f_{\mathbf{0}i}^{\mathbf{1}j} x_1^{\mathbf{0}i} + f_{\mathbf{1}i}^{\mathbf{1}j} x_1^{\mathbf{1}i} \quad (4.11)$$

The equation (4.9) follows from equations (4.2), (4.11). The equation (4.10) follows from equation (4.9) and comparing equations

$$(f_{\mathbf{0}i}^{\mathbf{0}j} + f_{\mathbf{0}i}^{\mathbf{1}j} i)(x_1^{\mathbf{0}i} + x_1^{\mathbf{1}i} i) = f_{\mathbf{0}i}^{\mathbf{0}j} x_1^{\mathbf{0}i} - f_{\mathbf{0}i}^{\mathbf{1}j} x_1^{\mathbf{1}i} + (f_{\mathbf{0}i}^{\mathbf{1}j} x_1^{\mathbf{0}i} + f_{\mathbf{0}i}^{\mathbf{0}j} x_1^{\mathbf{1}i}) i$$

$$\begin{pmatrix} f_{\mathbf{0}i}^{\mathbf{0}j} & -f_{\mathbf{0}i}^{\mathbf{1}j} \\ f_{\mathbf{0}i}^{\mathbf{1}j} & f_{\mathbf{0}i}^{\mathbf{0}j} \end{pmatrix} \begin{pmatrix} x_1^{\mathbf{0}i} \\ x_1^{\mathbf{1}i} \end{pmatrix} = \begin{pmatrix} f_{\mathbf{0}i}^{\mathbf{0}j} x_1^{\mathbf{0}i} - f_{\mathbf{0}i}^{\mathbf{1}j} x_1^{\mathbf{1}i} \\ f_{\mathbf{0}i}^{\mathbf{1}j} x_1^{\mathbf{0}i} + f_{\mathbf{0}i}^{\mathbf{0}j} x_1^{\mathbf{1}i} \end{pmatrix}$$

□

The theorem 4.2 describes the structure of $C$-linear mapping of $C$-algebra $A$. However, there exist $R$-linear mappings of $C$-algebra $A$, which are not $C$-linear mappings.

## 5. Antilinear Mapping of $C$-algebra

**Definition 5.1.** Morphism of $C$-vector spaces

$$I : z \in C \to z^* \in C \quad f : A_1 \to A_2 \quad (5.1)$$

is called antilinear mapping.[1] □

**Theorem 5.2.** *Let the mapping of $C$-vector spaces*

$$f : A_1 \to A_2$$

*be antilinear mapping. Then*

$$f(ca) = c^* f(a) \quad c \in C \quad a \in A \quad (5.2)$$

*Proof.* According to construction in section [4]-1.5, the mapping $I$ is morphism of $R$-algebra $C$. The equation (5.2) follows from the definition [5]-2.2.2, and the equation [5]-(2.2.4). □

If $c \in R$, then

$$f(ca) = cf(a)$$

because $c = c^*$. Therefore, antilinear mapping of $C$-algebra $A$ is $R$-linear mapping.

**Theorem 5.3** (The Cauchy-Riemann equations)**.** *Matrix of antilinear mapping of $C$-vector spaces*

$$f : A_1 \to A_2$$
$$y^{\mathbf{i}\mathbf{k}} = f_{\mathbf{j}\mathbf{m}}^{\mathbf{i}\mathbf{k}} x^{\mathbf{j}\mathbf{m}}$$

*satisfies relationship*

$$f_{\mathbf{0}i}^{\mathbf{0}j} = -f_{\mathbf{1}i}^{\mathbf{1}j}$$
$$f_{\mathbf{1}i}^{\mathbf{0}j} = f_{\mathbf{0}i}^{\mathbf{1}j} \quad (5.3)$$

[1] There is similar definition in [7], p. 234.

*Proof.* From equations (3.11), (3.14), (3.15) , it follows

$$(5.4)\quad \begin{aligned} (f(ca))^{0j} &= f^{0j}_{0i}(ca)^{0i} + f^{0j}_{1i}(ca)^{1i} \\ &= f^{0j}_{0i}(c^0a^{0i} - c^1a^{1i}) + f^{0j}_{1i}(c^0a^{1i} + c^1a^{0i}) \\ &= f^{0j}_{0i}c^0a^{0i} - f^{0j}_{0i}c^1a^{1i} + f^{0j}_{1i}c^0a^{1i} + f^{0j}_{1i}c^1a^{0i} \end{aligned}$$

$$(5.5)\quad \begin{aligned} (f(ca))^{1j} &= f^{1j}_{0i}(ca)^{0i} + f^{1j}_{1i}(ca)^{1i} \\ &= f^{1j}_{0i}(c^0a^{0i} - c^1a^{1i}) + f^{1j}_{1i}(c^0a^{1i} + c^1a^{0i}) \\ &= f^{1j}_{0i}c^0a^{0i} - f^{1j}_{0i}c^1a^{1i} + f^{1j}_{1i}c^0a^{1i} + f^{1j}_{1i}c^1a^{0i} \end{aligned}$$

From equations (5.2), (3.11), (3.13), (3.14), (3.15) , it follows

$$(5.6)\quad \begin{aligned} (c^*f(a))^{0j} &= c^0(f(a))^{0j} + c^1(f(a))^{1j} \\ &= c^0(f^{0j}_{0i}a^{0i} + f^{0j}_{1i}a^{1i}) + c^1(f^{1j}_{0i}a^{0i} + f^{1j}_{1i}a^{1i}) \\ &= c^0f^{0j}_{0i}a^{0i} + c^0f^{0j}_{1i}a^{1i} + c^1f^{1j}_{0i}a^{0i} + c^1f^{1j}_{1i}a^{1i} \end{aligned}$$

$$(5.7)\quad \begin{aligned} (c^*f(a))^{1j} &= c^0(f(a))^{1j} - c^1(f(a))^{0j} \\ &= c^0(f^{1j}_{0i}a^{0i} + f^{1j}_{1i}a^{1i}) - c^1(f^{0j}_{0i}a^{0i} + f^{0j}_{1i}a^{1i}) \\ &= c^0f^{1j}_{0i}a^{0i} + c^0f^{1j}_{1i}a^{1i} - c^1f^{0j}_{0i}a^{0i} - c^1f^{0j}_{1i}a^{1i} \end{aligned}$$

From equations (5.2), (5.4), (5.6) , it follows

$$(5.8)\quad \begin{aligned} &f^{0j}_{0i}c^0a^{0i} - f^{0j}_{0i}c^1a^{1i} + f^{0j}_{1i}c^0a^{1i} + f^{0j}_{1i}c^1a^{0i} \\ =&c^0f^{0j}_{0i}a^{0i} + c^0f^{0j}_{1i}a^{1i} + c^1f^{1j}_{0i}a^{0i} + c^1f^{1j}_{1i}a^{1i} \end{aligned}$$

From equations (5.2), (5.5), (5.7) , it follows

$$(5.9)\quad \begin{aligned} &f^{1j}_{0i}c^0a^{0i} - f^{1j}_{0i}c^1a^{1i} + f^{1j}_{1i}c^0a^{1i} + f^{1j}_{1i}c^1a^{0i} \\ =&c^0f^{1j}_{0i}a^{0i} + c^0f^{1j}_{1i}a^{1i} - c^1f^{0j}_{0i}a^{0i} - c^1f^{0j}_{1i}a^{1i} \end{aligned}$$

Since $c^k$, $a^{ki}$ are arbitrary, then the equation (5.3) follows from equations (5.8), (5.9). □

**Theorem 5.4.** *Let $\overline{\overline{e}}_{iC}$ be $C$-basis of $C$-vector space $A_i$, $i = 1, 2$. Let $\overline{\overline{e}}_{iR}$ be corresponding $R$-basis of $C$-vector space $A_i$. The matrix of antilinear mapping*

$$f : A_1 \to A_2$$

*has form*

$$(5.10)\quad \begin{aligned} x_2^{0j} &= f^{0j}_{0i}x_1^{0i} + f^{1j}_{0i}x_1^{1i} \\ x_2^{1j} &= f^{1j}_{0i}x_1^{0i} - f^{0j}_{0i}x_1^{1i} \end{aligned}$$

*relative to bases $\overline{\overline{e}}_{1R}$, $\overline{\overline{e}}_{2R}$ and has form*

$$(5.11)\quad x_2^j = f^j_i I \circ x_1^i \quad I \circ z = z^* \quad f^j_i = f^{0j}_{0i} + f^{1j}_{0i}i$$

*relative to bases $\overline{\overline{e}}_{1C}$, $\overline{\overline{e}}_{2C}$*

*Proof.* Since the mapping $f$ is $R$-linear mapping, then relative to bases $\overline{\overline{e}}_{1R}$, $\overline{\overline{e}}_{2R}$ mapping $f$ has form

$$
\begin{aligned}
x_2^{0j} &= f_{0i}^{0j} x_1^{0i} + f_{1i}^{0j} x_1^{1i} \\
x_2^{1j} &= f_{0i}^{1j} x_1^{0i} + f_{1i}^{1j} x_1^{1i}
\end{aligned}
\tag{5.12}
$$

The equation (5.10) follows from equations (5.3), (5.12). The equation (5.11) follows from the equation (5.10) and comparing of equations

$$
\begin{aligned}
(f_{0i}^{0j} + f_{0i}^{1j} i) I \circ (x_1^{0i} + x_1^{1i} i) &= (f_{0i}^{0j} + f_{0i}^{1j} i)(x_1^{0i} - x_1^{1i} i) \\
&= f_{0i}^{0j} x_1^{0i} + f_{0i}^{1j} x_1^{1i} + (f_{0i}^{1j} x_1^{0i} - f_{0i}^{0j} x_1^{1i}) i
\end{aligned}
$$

$$
\begin{pmatrix} f_{0i}^{0j} & f_{0i}^{1j} \\ f_{0i}^{1j} & -f_{0i}^{0j} \end{pmatrix}
\begin{pmatrix} x_1^{0i} \\ x_1^{1i} \end{pmatrix}
=
\begin{pmatrix} f_{0i}^{0j} x_1^{0i} + f_{0i}^{1j} x_1^{1i} \\ f_{0i}^{1j} x_1^{0i} - f_{0i}^{0j} x_1^{1i} \end{pmatrix}
$$

□

## 6. Morphism of Algebra

In this section, we consider algebra $A$ over ring $D$. According to construction that was done in subsections [5]-4.4.2, [5]-4.4.3, a diagram of representations of $D$-algebra has form

$$
\begin{array}{ccccc}
D & \xrightarrow[*]{f_{1,2}} & A & \xrightarrow[*]{f_{2,3}} & A \\
 & & & & \uparrow {\scriptstyle *}\, f_{1,2} \\
 & & & & D
\end{array}
\qquad
\begin{aligned}
f_{1,2}(d) &: \overline{v} \rightarrow d\,\overline{v} \\
f_{2,3}(\overline{v}) &: \overline{w} \rightarrow \overline{C}(\overline{v}, \overline{w}) \\
&\overline{C} \in \mathcal{L}(A^2; A)
\end{aligned}
\tag{6.1}
$$

On the diagram of representations (6.1), $D$ is ring, $A$ is Abelian group. We initially consider the vertical representation, and then we consider the horizontal representation. We assume that $A$ is free $D\star$-module ([1], p. 135). Therefore, vectors of basis $\overline{\overline{e}}$ of $D\star$-module $A$ are linear independent. Bilinear mapping

$$\overline{C} : A \times A \rightarrow A$$

has following representation relative basis $\overline{\overline{e}}$

$$\overline{C}(\overline{v}, \overline{w}) = C_{ij}^{k} v^{i} w^{j} \overline{e}_{k} \tag{6.2}$$

According to the theorem [5]-4.4.1, bilinear mapping $\overline{C}$ generates product in $D\star$-algebra $A$

$$\overline{C}(\overline{v}, \overline{w}) = \overline{v}\,\overline{w} \tag{6.3}$$

and coordinates of bilinear mapping $\overline{C}$ relative to the basis $\overline{\overline{e}}$ are structural constants of this product.

**Theorem 6.1.** *Let diagram of representations*

$$
\begin{array}{ccccc}
D_1 & \xrightarrow[*]{f_{1\cdot 1,2}} & A_1 & \xrightarrow[*]{f_{1\cdot 2,3}} & A_1 \\
 & & & & \uparrow {\scriptstyle *}\, f_{1\cdot 1,2} \\
 & & & & D_1
\end{array}
\qquad
\begin{aligned}
f_{1\cdot 1,2}(d) &: \overline{v} \rightarrow d\,\overline{v} \\
f_{1\cdot 2,3}(\overline{v}) &: \overline{w} \rightarrow \overline{C}_1(\overline{v}, \overline{w}) \\
&\overline{C}_1 \in \mathcal{L}(A_1^2; A_1)
\end{aligned}
\tag{6.4}
$$

*describe $D_1\star$-algebra $A_1$. Let diagram of representations*

$$
\begin{array}{ccccc}
D_2 & \xrightarrow[*]{f_{2\cdot 1,2}} & A_2 & \xrightarrow[*]{f_{2\cdot 2,3}} & A_2 \\
 & & & & \uparrow{\scriptstyle *}\, f_{2\cdot 1,2} \\
 & & & & D_2
\end{array}
\qquad
\begin{array}{l}
f_{2\cdot 1,2}(d):\overline{v}\rightarrow d\,\overline{v}\\
f_{2\cdot 2,3}(\overline{v}):\overline{w}\rightarrow \overline{C}_2(\overline{v},\overline{w})\\
\overline{C}_2\in\mathcal{L}(A_2^2;A_2)
\end{array}
\tag{6.5}
$$

*describe $D_2\star$-algebra $A_2$. Morphism*[2] *of $D_1\star$-algebra $A_1$ into $D_2\star$-algebra $A_2$*

$$r_1 : D_1 \rightarrow D_2 \quad \overline{r}_2 : A_1 \rightarrow A_2$$

*is $D\star$-linear mapping of $D_1\star$-algebra $A_1$ into $D_2\star$-algebra $A_2$ such that*

$$\overline{r}_2(\overline{a}\overline{b}) = \overline{r}_2(\overline{a})\overline{r}_2(\overline{b}) \tag{6.6}$$

*Proof.* According to equation [5]-(4.2.3), morphism $(r_1, \overline{r}_2)$ of representation $f_{1,2}$ satisfies to the equation

$$
\begin{aligned}
\overline{r}_2(f_{1\cdot 1,2}(d)(\overline{a})) &= f_{2\cdot 1,2}(r_1(d))(\overline{r}_2(\overline{a}))\\
\overline{r}_2(d\,\overline{a}) &= r_1(d)\overline{r}_2(\overline{a})
\end{aligned}
\tag{6.7}
$$

Therefore, the mapping $(r_1, \overline{r}_2)$ is $D\star$-linear mapping.

According to equation [5]-(4.2.3), the morphism $(\overline{r}_2, \overline{r}_2)$ of representation $f_{2,3}$ satisfies to the equation[3]

$$\overline{r}_2(f_{1\cdot 2,3}(\overline{a}_2)(\overline{a}_3)) = f_{2\cdot 2,3}(\overline{r}_2(\overline{a}_2))(\overline{r}_2(\overline{a}_3)) \tag{6.8}$$

From equations (6.8), (6.4), (6.5), it follows that

$$\overline{r}_2(\overline{C}_1(\overline{v},\overline{w})) = \overline{C}_2(\overline{r}_2(\overline{v}), \overline{r}_2(\overline{w})) \tag{6.9}$$

Equation (6.6) follows from equations (6.9), (6.3). □

**Definition 6.2.** The morphism of representations of $D_1\star$-algebra $A_1$ into $D_2\star$-algebra $A_2$ is called **$D\star$-linear homomorphism** of $D_1\star$-algebra $A_1$ into $D_2\star$-algebra $A_2$. □

**Theorem 6.3.** *Let $\overline{\overline{e}}_1$ be the basis of $D_1\star$-algebra $A_1$. Let $\overline{\overline{e}}_2$ be the basis of $D_2\star$-algebra $A_2$. Then linear homomorphism*[4] *$(r_1, \overline{r}_2)$ of $D_1\star$-algebra $A_1$ into $D_2\star$-algebra $A_2$ has presentation*

$$
\begin{aligned}
\overline{b} &= r_1(a)^*{}_*r_2{}^*{}_*\overline{e}_2 = r_1(a^j)r_{2\cdot}{}^i_j\,\overline{e}_{2\cdot i}\\
b &= r_1(a)^*{}_*r_2
\end{aligned}
\tag{6.10}
$$

*relative to selected bases. Here*

- *$a$ is coordinate matrix of vector $\overline{a}$ relative the basis $\overline{\overline{e}}_1$.*

---

[2] In the remark [2]-4.4.5, I show that, when we study linear mapping, the choice of mapping $r_1$ may be limited to mapping

$$r_1 : D \rightarrow D \quad r_1(d) = d$$

However, in the section 10, I consider antilinear mapping of the algebra. That is why I wrote theorems in this section so that we can compare these theorems with theorems in section 10.

[3] Since in diagrams of representations (6.4), (6.5), supports of $\Omega_2$-algebra and $\Omega_3$-algebra coincide, then morphisms of representations on levels 2 and 3 coincide also.

[4] This theorem is similar to the theorem [2]-4.4.3.

- *$b$ is coordinate matrix of vector*
$$\overline{b}=\overline{r}_2(\overline{a})$$
*relative the basis $\overline{\overline{e}}_2$.*
- *$r_2$ is coordinate matrix of set of vectors $(\overline{r}_2(\overline{e}_{1\cdot i}))$ relative the basis $\overline{\overline{e}}_2$. The matrix $r_2$ is called* **matrix of linear homomorphism** *relative bases $\overline{\overline{e}}_1$ and $\overline{\overline{e}}_2$.*

*Proof.* Vector $\overline{a}\in A_1$ has expansion
$$\overline{a}=a^*{}_*\overline{e}_1$$
relative to the basis $\overline{\overline{e}}_1$. Vector $\overline{b}\in A_2$ has expansion
$$\overline{b}=b^*{}_*\overline{e}_2 \tag{6.11}$$
relative to the basis $\overline{\overline{e}}_2$.

Since $(r_1,\overline{r}_2)$ is a linear homomorphism, then from (6.7), it follows that
$$\overline{b}=\overline{r}_2(\overline{a})=\overline{r}_2(a^*{}_*\overline{e}_1)=r_1(a)^*{}_*\overline{r}_2(\overline{e}_1) \tag{6.12}$$
where
$$r_1(a)=\begin{pmatrix}r_1(a^1)\\ \dots\\ r_1(a^n)\end{pmatrix}$$
$\overline{r}_2(\overline{e}_{1\cdot i})$ is also a vector of $D\star$-module $A_2$ and has expansion
$$\overline{r}_2(\overline{e}_{1\cdot i})=r_{2\cdot i}{}^*{}_*\overline{e}_2=r_{2\cdot}{}^j_i\,\overline{e}_j \tag{6.13}$$
relative to basis $\overline{\overline{e}}_2$. Combining (6.12) and (6.13) we get
$$\overline{b}=r_1(a)^*{}_*r_2{}^*{}_*\overline{e}_2 \tag{6.14}$$
(6.10) follows from comparison of (6.11) and (6.14) and theorem [2]-4.3.3. □

**Theorem 6.4.** *Let $\overline{\overline{e}}_1$ be the basis of $D_1\star$-algebra $A_1$. Let $\overline{\overline{e}}_2$ be the basis of $D_2\star$-algebra $A_2$. If the mapping $r_1$ is injection, then there is relation between the matrix of linear homomorphism and structural constants*
$$r_1(C_{1\cdot}{}^k_{ij})r_{2\cdot}{}^l_k=r_{2\cdot}{}^p_i r_{2\cdot}{}^q_j C_{2\cdot}{}^l_{pq} \tag{6.15}$$

*Proof.* Let
$$\overline{a},\overline{b}\in A_1\quad \overline{a}=a^*{}_*\overline{e}_1\quad \overline{b}=b^*{}_*\overline{e}_1$$
From equations (6.2), (6.3), (6.4), it follows that
$$\overline{a}\overline{b}=a^ib^jC_{1\cdot}{}^k_{ij}\overline{e}_{1\cdot k} \tag{6.16}$$
From equations (6.7), (6.16), it follows that
$$\overline{r}_2(\overline{a}\overline{b})=r_1(C_{1\cdot}{}^k_{ij}a^ib^j)\overline{r}_2(\overline{e}_{1\cdot k}) \tag{6.17}$$
Since the mapping $r_1$ is homomorphism of rings, then from equation (6.17), it follows that
$$\overline{r}_2(\overline{a}\overline{b})=r_1(C_{1\cdot}{}^k_{ij})r_1(a^i)r_1(b^j)\overline{r}_2(\overline{e}_{1\cdot k}) \tag{6.18}$$
From the theorem 6.3 and the equation (6.18), it follows that
$$\overline{r}_2(\overline{a}\overline{b})=r_1(C_{1\cdot}{}^k_{ij})r_1(a^i)r_1(b^j)r_{2\cdot}{}^l_k\,\overline{e}_{2\cdot l} \tag{6.19}$$

From the equation (6.6) and the theorem 6.3, it follows that

$$\overline{r}_2(\overline{a}\overline{b}) = \overline{r}_2(\overline{a})\overline{r}_2(\overline{b}) = r_1(a^i)r_{2\cdot}{}^p_i\,\overline{e}_{2\cdot p}r_1(b^j)r_{2\cdot}{}^q_j\,\overline{e}_{2\cdot q} \tag{6.20}$$

From equations (6.2), (6.3), (6.5), (6.20), it follows that

$$\overline{r}_2(\overline{a}\overline{b}) = r_1(a^i)r_{2\cdot}{}^p_i r_1(b^j)r_{2\cdot}{}^q_j C_{2\cdot}{}^l_{pq}\,\overline{e}_{2\cdot l} \tag{6.21}$$

From equations (6.19), (6.21), it follows that

$$r_1(C_{1\cdot}{}^k_{ij})r_1(a^i)r_1(b^j)r_{2\cdot}{}^l_k\,\overline{e}_{2\cdot l} = r_1(a^i)r_{2\cdot}{}^p_i r_1(b^j)r_{2\cdot}{}^q_j C_{2\cdot}{}^l_{pq}\,\overline{e}_{2\cdot l} \tag{6.22}$$

The equation (6.15) follows from the equation (6.22), because vectors of basis $\overline{\overline{e}}_2$ are linear independent, and $a^i$, $b^i$ (and therefore, $r_1(a^i)$, $r_1(b^i)$) are arbitrary values. □

## 7. Linear Automorphism of Quaternion Algebra

Determining of coordinates of linear automorphism is not a simple task. In this section we consider example of nontrivial linear automorphism of quaternion algebra.

**Theorem 7.1.** *Coordinates of linear automorphism of quaternion algebra satisfy to the system of equations*

$$\begin{cases} r^1_1 = r^2_2r^3_3 - r^3_2r^2_3 & r^1_2 = r^2_3r^3_1 - r^3_3r^2_1 & r^1_3 = r^2_1r^3_2 - r^3_1r^2_2 \\ r^2_1 = r^3_2r^1_3 - r^1_2r^3_3 & r^2_2 = r^3_3r^1_1 - r^1_3r^3_1 & r^2_3 = r^3_1r^1_2 - r^1_1r^3_2 \\ r^3_1 = r^1_2r^2_3 - r^2_2r^1_3 & r^3_2 = r^1_3r^2_1 - r^2_3r^1_1 & r^3_3 = r^1_1r^2_2 - r^2_1r^1_2 \end{cases} \tag{7.1}$$

*Proof.* According to the theorems [4]-3.3.1, 6.4, linear automorphism of quaternion algebra satisfies to equations

$$\begin{array}{llll} r^l_0 = r^p_0r^q_0C^l_{pq} & r^l_1 = r^p_0r^q_1C^l_{pq} & r^l_2 = r^p_0r^q_2C^l_{pq} & r^l_3 = r^p_0r^q_3C^l_{pq} \\ r^l_1 = r^p_1r^q_0C^l_{pq} & -r^l_0 = r^p_1r^q_1C^l_{pq} & r^l_3 = r^p_1r^q_2C^l_{pq} & -r^l_2 = r^p_1r^q_3C^l_{pq} \\ r^l_2 = r^p_2r^q_0C^l_{pq} & -r^l_3 = r^p_2r^q_1C^l_{pq} & -r^l_0 = r^p_2r^q_2C^l_{pq} & r^l_0 = r^p_2r^q_3C^l_{pq} \\ r^l_3 = r^p_3r^q_0C^l_{pq} & r^l_2 = r^p_3r^q_1C^l_{pq} & -r^l_1 = r^p_3r^q_2C^l_{pq} & -r^l_0 = r^p_3r^q_3C^l_{pq} \end{array} \tag{7.2}$$

From the equation (7.2), it follows that

$$\begin{array}{lll} r^l_1 = r^p_0r^q_1C^l_{pq} = r^p_2r^q_3C^l_{pq} & r^p_0r^q_1C^l_{pq} = r^p_0r^q_1C^l_{qp} & r^p_2r^q_3C^l_{pq} = -r^p_2r^q_3C^l_{qp} \\ r^l_2 = r^p_0r^q_2C^l_{pq} = r^p_3r^q_1C^l_{pq} & r^p_0r^q_2C^l_{pq} = r^p_0r^q_2C^l_{qp} & r^p_3r^q_1C^l_{pq} = -r^p_1r^q_3C^l_{pq} \\ r^l_3 = r^p_0r^q_3C^l_{pq} = r^p_1r^q_2C^l_{pq} & r^p_0r^q_3C^l_{pq} = r^p_0r^q_3C^l_{qp} & r^p_1r^q_2C^l_{pq} = -r^p_1r^q_2C^l_{qp} \end{array} \tag{7.3}$$

$$r^l_0 = r^p_0r^q_0C^l_{pq} = -r^p_1r^q_1C^l_{pq} = -r^p_2r^q_2C^l_{pq} = -r^p_3r^q_3C^l_{pq} \tag{7.4}$$

If $l = 0$, then from the equation

$$C^0_{pq} = C^0_{qp}$$

it follows that

$$r^p_ir^q_jC^0_{pq} = r^p_ir^q_jC^0_{qp} \tag{7.5}$$

From the equation (7.3) for $l=0$ and the equation (7.5), it follows that

$$r_1^0=r_2^0=r_3^0=0 \tag{7.6}$$

If $l=1, 2, 3$, then we can write the equation (7.3) in the following form

$$\left\{\begin{array}{l}
r_i^l=r_0^l r_i^0 C_{l0}^l+r_0^0 r_i^l C_{0l}^l+r_0^a r_i^b C_{ab}^l+r_0^b r_i^a C_{ba}^l\\
\quad r_0^l r_i^0 C_{l0}^l+r_0^0 r_i^l C_{0l}^l+r_0^a r_i^b C_{ab}^l+r_0^b r_i^a C_{ba}^l\\
\quad =r_0^l r_i^0 C_{0l}^l+r_0^0 r_i^l C_{l0}^l+r_0^a r_i^b C_{ba}^l+r_0^b r_i^a C_{ab}^l\\
\quad i=1,2,3\\
r_i^l=r_k^0 r_j^l C_{0l}^l+r_k^l r_j^0 C_{l0}^l+r_k^a r_j^b C_{ab}^l+r_k^b r_j^a C_{ba}^l\\
\quad r_k^0 r_j^l C_{0l}^l+r_k^l r_j^0 C_{l0}^l+r_k^a r_j^b C_{ab}^l+r_k^b r_j^a C_{ba}^l\\
\quad =-r_k^0 r_j^l C_{l0}^l-r_k^l r_j^0 C_{0l}^l-r_k^a r_j^b C_{ba}^l-r_k^b r_j^a C_{ab}^l\\
\qquad i=1 \quad k=2 \quad j=3\\
\qquad i=2 \quad k=3 \quad j=1\\
\qquad i=3 \quad k=1 \quad j=2\\
\quad 0<a<b \quad a\ne l \quad b\ne l
\end{array}\right. \tag{7.7}$$

From equations (7.7), (7.6) and equations

$$\begin{array}{c}C_{0l}^l=C_{l0}^l=1\\ C_{ab}^l=-C_{ba}^l\end{array} \tag{7.8}$$

it follows that

$$\left\{\begin{array}{l}
r_i^l=r_0^0 r_i^l+r_0^a r_i^b C_{ab}^l-r_0^b r_i^a C_{ab}^l\\
\quad r_0^0 r_i^l+r_0^a r_i^b C_{ab}^l-r_0^b r_i^a C_{ab}^l\\
\quad =r_0^0 r_i^l-r_0^a r_i^b C_{ab}^l+r_0^b r_i^a C_{ab}^l\\
\quad i=1,2,3\\
r_i^l=r_k^a r_j^b C_{ab}^l-r_k^b r_j^a C_{ab}^l\\
\quad r_k^a r_j^b C_{ab}^l-r_k^b r_j^a C_{ab}^l\\
\quad =r_k^a r_j^b C_{ab}^l-r_k^b r_j^a C_{ab}^l\\
\qquad i=1 \quad k=2 \quad j=3\\
\qquad i=2 \quad k=3 \quad j=1\\
\qquad i=3 \quad k=1 \quad j=2\\
\quad 0<a<b \quad a\ne l \quad b\ne l
\end{array}\right. \tag{7.9}$$

From equations (7.9) it follows that

$$(7.10)\quad \left\{\begin{array}{l} r_i^l = r_0^0 r_i^l \\ \qquad r_0^a r_i^b - r_0^b r_i^a = 0 \\ \qquad i = 1, 2, 3 \\ r_i^l = r_k^a r_j^b C_{ab}^l - r_k^b r_j^a C_{ab}^l \\ \qquad i = 1 \quad k = 2 \quad j = 3 \\ \qquad i = 2 \quad k = 3 \quad j = 1 \\ \qquad i = 3 \quad k = 1 \quad j = 2 \\ \qquad 0 < a < b \quad a \neq l \quad b \neq l \end{array}\right.$$

From the equation (7.10) it follows that

$$(7.11)\quad r_0^0 = 1$$

From the equation (7.4) for $l = 0$, it follows that

$$(7.12)\quad \begin{aligned} r_0^0 &= r_0^0 r_0^0 - r_0^1 r_0^1 - r_0^2 r_0^2 - r_0^3 r_0^3 \\ &= -r_i^0 r_i^0 + r_i^1 r_i^1 + r_i^2 r_i^2 + r_i^3 r_i^3 \\ i &= 1, 2, 3 \end{aligned}$$

From equations (7.6), (7.10), (7.12), it follows that

$$(7.13)\quad \begin{aligned} 0 &= r_0^1 r_0^1 + r_0^2 r_0^2 + r_0^3 r_0^3 \\ 1 &= r_i^1 r_i^1 + r_i^2 r_i^2 + r_i^3 r_i^3 \\ i &= 1, 2, 3 \end{aligned}$$

From equations (7.13) it follows that[5]

$$(7.14)\quad r_0^1 = r_0^2 = r_0^3 = 0$$

From the equation (7.4) for $l > 0$, it follows that

$$(7.15)\quad \begin{aligned} r_0^l &= r_0^l r_0^0 C_{l0}^l + r_0^0 r_0^l C_{0l}^l + r_0^a r_0^b C_{ab}^l + r_0^b r_0^a C_{ba}^l \\ &= -r_i^l r_i^0 C_{l0}^l - r_i^0 r_i^l C_{0l}^l - r_i^a r_i^b C_{ab}^l - r_i^b r_i^a C_{ba}^l \\ & i > 0 \\ & l > 0 \quad 0 < a < b \quad a \neq l \quad b \neq l \end{aligned}$$

[5]Here, we rely on the fact that quaternion algebra is defined over real field. If we consider quaternion algebra over complex field, then the equation (7.13) defines a cone in the complex space. Correspondingly, we have wider choice of coordinates of linear automorphism.

Equations (7.15) are identically true by equations (7.6), (7.14), (7.8). From equations (7.14), (7.10), it follows that

$$(7.16)\quad \begin{cases} r_i^l = r_k^a r_j^b C_{ab}^l - r_k^b r_j^a C_{ab}^l \\ \quad i=1 \quad k=2 \quad j=3 \\ \quad i=2 \quad k=3 \quad j=1 \\ \quad i=3 \quad k=1 \quad j=2 \\ \quad l>0 \quad 0<a<b \quad a\neq l \quad b\neq l \end{cases}$$

Equations (7.1) follow from equations (7.16). □

**Example 7.2.** It is evident that coordinates

$$r_j^i = \delta_j^i$$

satisfy the equation (7.1). It can be verified directly that coordinates of mapping

$$r_0^0 = 1 \quad r_2^1 = 1 \quad r_3^2 = 1 \quad r_1^3 = 1$$

also satisfy the equation (7.1). The matrix of coordinates of this mapping has form

$$r = \begin{pmatrix} 1 & 0 & 0 & \\ 0 & 0 & 1 & 0 \\ 0 & 0 & 0 & 1 \\ 0 & 1 & 0 & 0 \end{pmatrix}$$

According to the theorem [4]-3.3.4, standard components of the mapping $r$ have form

$$\begin{matrix} r^{00} = \frac{1}{4} & r^{11} = -\frac{1}{4} & r^{22} = -\frac{1}{4} & r^{33} = -\frac{1}{4} \\ r^{10} = -\frac{1}{4} & r^{01} = \frac{1}{4} & r^{32} = -\frac{1}{4} & r^{23} = -\frac{1}{4} \\ r^{20} = -\frac{1}{4} & r^{31} = -\frac{1}{4} & r^{02} = \frac{1}{4} & r^{13} = -\frac{1}{4} \\ r^{30} = -\frac{1}{4} & r^{21} = -\frac{1}{4} & r^{12} = -\frac{1}{4} & r^{03} = \frac{1}{4} \end{matrix}$$

Therefore, the mapping $r$ has form

$$r(a) = a^0 + a^2 i + a^3 j + a^1 k$$

$$\begin{aligned} r(a) = \frac{1}{4}(&a - iai - jaj - kak - ia + ai - kaj - jak \\ &- ja - kai + aj - iak - ka - jai - iaj + ak) \end{aligned}$$

□

## 8. Ring with Conjugation

Let there exist commutative subring $F$ of ring $D$ such that $F \neq D$ and ring $D$ is a free finite dimensional module over the ring $F$. Let $\overline{\overline{e}}$ be the basis of free module $D$ over ring $F$. We assume that $\overline{e}_0 = 1$.

**Theorem 8.1.** *The ring $D$ is $F$-algebra.*

*Proof.* Consider mapping

$$f : D \times D \to D \quad a, b \in D \to ab \in D \tag{8.1}$$

Since the product is distributive over addition, then

$$\begin{aligned} f(a_1 + a_2, b) &= (a_1 + a_2)b = a_1 b + a_2 b = f(a_1, b) + f(a_2, b) \\ f(a, b_1 + b_2) &= a(b_1 + b_2) = ab_1 + ab_2 = f(a, b_1) + f(a, b_2) \end{aligned} \tag{8.2}$$

If $a \in F$, then

$$\begin{aligned} f(ab, c) &= (ab)c = a(bc) = af(b, c) \\ f(b, ac) &= b(ac) = (ba)c = (ab)c = a(bc) = af(b, c) \end{aligned} \tag{8.3}$$

From equations (8.2), (8.3), it follows that the mapping $f$ is bilinear mapping over the ring $F$. According to the definition [4]-2.2.1, the ring $D$ is $F$-algebra. □

**Corollary 8.2.** *The ring $F$ is subalgebra of $F$-algebra $D$.* □

Consider mappings

$$\begin{aligned} \mathrm{Re} &: D \to D \\ \mathrm{Im} &: D \to D \end{aligned}$$

defined by equation

$$\mathrm{Re}\, d = d^0 \quad \mathrm{Im}\, d = d - d^0 \quad d \in D \quad d = d^i \overline{e}_i \tag{8.4}$$

The expression $\mathrm{Re}\, d$ is called **scalar of element** $d$. The expression $\mathrm{Im}\, d$ is called **vector of element** $d$.

According to (8.4)

$$F = \{d \in D : \mathrm{Re}\, d = d\}$$

We will use notation $\mathrm{Re}\, D$ to denote **scalar algebra of ring** $D$.

**Theorem 8.3.** *The set*

$$\mathrm{Im}\, D = \{d \in D : \mathrm{Re}\, d = 0\} \tag{8.5}$$

*is* $(\mathrm{Re}\, D)$*-module which is called* **vector module of ring** $D$.

$$D = \mathrm{Re}\, D \oplus \mathrm{Im}\, D \tag{8.6}$$

*Proof.* Let $c$, $d \in \mathrm{Im}\, D$. Then $c_0 = d_0 = 0$. Therefore,

$$(c + d)_0 = c_0 + d_0 = 0$$

If $a \in \mathrm{Re}\, D$, then

$$(ad)_0 = ad_0 = 0$$

Therefore, $\mathrm{Im}\, D$ is $(\mathrm{Re}\, D)$-module.

Sequence of modules

$$0 \longrightarrow \mathrm{Re}\, D \xrightarrow{\mathrm{id}} D \xrightarrow{\mathrm{Im}} \mathrm{Im}\, D \longrightarrow 0$$

is exact sequence. According to the definition (8.4) of the mapping Re , following diagram is commutative

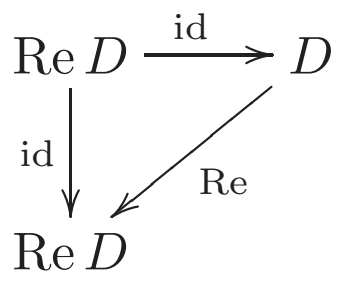

According to the statement 2 of the proposition [1]-III.3.2,

$$D = \mathrm{id}(\mathrm{Re}\,D) \oplus \ker \mathrm{Re} \tag{8.7}$$

According to the definition (8.5)

$$\ker \mathrm{Re} = \{d \in D : \mathrm{Re}\,d = 0\} = \mathrm{Im}\,D \tag{8.8}$$

The equation (8.6) follows from the equations (8.7), (8.8). □

According to the theorem 8.3, there is unique defined representation

$$d = \mathrm{Re}\,d + \mathrm{Im}\,d \tag{8.9}$$

**Definition 8.4.** The mapping

$$d^* = \mathrm{Re}\,d - \mathrm{Im}\,d \tag{8.10}$$

is called **conjugation in ring** provided that this mapping satisfies

$$(cd)^* = d^*\,c^* \tag{8.11}$$

The ring $D$ equipped with conjugation is called **ring with conjugation**. □

According to the theorem 8.1, the ring with conjugation $D$ is associative $(\mathrm{Re}\,D)$-algebra, which we also call **algebra with conjugation**.

**Corollary 8.5.**

$$(d^*)^* = d \tag{8.12}$$

□

It is evident that there exist rings where condition (8.11) is not fulfilled. To understand the properties of a ring with conjugation we consider structural constants of $(\mathrm{Re}\,D)$-algebra $D$. The equation

$$C_{0k}^{l} = C_{k0}^{l} = \delta_{l}^{k} \tag{8.13}$$

follows from the equation $1d = d1 = d$.

**Theorem 8.6.** *The ring $D$ is ring with conjugation if*

$$C_{kl}^{0} = C_{lk}^{0} \quad C_{kl}^{p} = -C_{lk}^{p} \tag{8.14}$$

$$\boldsymbol{1 < k < n \quad 1 < l < n \quad 1 < p < n}$$

*Proof.* From equations (8.9), (8.10) it follows that

$$\begin{aligned}(cd)^* &= (\mathrm{Re}\,c\ \mathrm{Re}\,d + \mathrm{Re}\,c\ \mathrm{Im}\,d + \mathrm{Im}\,c\ \mathrm{Re}\,d + \mathrm{Im}\,c\ \mathrm{Im}\,d)^* \\ &= \mathrm{Re}\,c\ \mathrm{Re}\,d - \mathrm{Re}\,c\ \mathrm{Im}\,d - \mathrm{Im}\,c\ \mathrm{Re}\,d + (\mathrm{Im}\,c\ \mathrm{Im}\,d)^*\end{aligned} \tag{8.15}$$

$$d^*\,c^* = (\mathrm{Re}\,d - \mathrm{Im}\,d)(\mathrm{Re}\,c - \mathrm{Im}\,c)$$
$$= \mathrm{Re}\,d\ \mathrm{Re}\,c - \mathrm{Re}\,d\ \mathrm{Im}\,c - \mathrm{Im}\,d\ \mathrm{Re}\,c + \mathrm{Im}\,d\ \mathrm{Im}\,c \tag{8.16}$$

Therefore, from the equation (8.11), it follows that

$$(\mathrm{Im}\,c\ \mathrm{Im}\,d)^* = \mathrm{Im}\,d\ \mathrm{Im}\,c \tag{8.17}$$

Let $\overline{\overline{e}}$ be the basis of the $(\mathrm{Re}\,D)$-algebra $D$. From the equation (8.17), it follows that

$$(C^0_{kl}c^k d^l + C^p_{kl}c^k d^l \overline{e}_p)^* = C^0_{kl}c^k d^l - C^p_{kl}c^k d^l \overline{e}_p$$
$$= C^0_{kl}d^k c^l + C^p_{kl}d^k c^l \overline{e}_p \tag{8.18}$$

$$1 < k < n \quad 1 < l < n \quad 1 < p < n$$

The equation (8.14) follows from the equation (8.18). □

**Corollary 8.7.**

$$\overline{e}_k\,\overline{e}_k \in \mathrm{Re}\,D$$

□

We can represent the conjugation using the matrix $I$

$$d^* = I\circ d \quad I^0_k = \delta^0_k \quad I^k_0 = \delta^k_0 \quad k = 0, ..., n$$
$$I^m_k = -\delta^m_k \quad I^k_m = -\delta^k_m \quad m = 1, ..., n \tag{8.19}$$

**Example 8.8.** The product in the set of complex numbers $C$ is commutative. However, complex field has subfield $R$. The vector space $\mathrm{Im}\,C$ has dimension 1 and the basis $\overline{e}_1 = i$. Accordingly

$$c^* = c^0 - c^1 i$$

□

**Example 8.9.** The division ring of quaternions $H$ has subfield $R$. The vector space $\mathrm{Im}\,H$ has dimension 3 and the basis

$$\overline{e}_1 = i \quad \overline{e}_2 = j \quad \overline{e}_3 = k$$

Accordingly

$$d^* = d^0 - d^1 i - d^2 j - d^3 k$$

□

## 9. Linear Mapping of Algebra with Conjugation

We can consider the ring $D$ as $D\star$-module of dimension 1.

**Theorem 9.1.** *$D\star$-linear mapping of the ring with conjugation $D$ has form*

$$f(c) = cd \quad d \in D \tag{9.1}$$

*Let $\overline{\overline{e}}$ be $(\mathrm{Re}\,D)$-basis of the ring $D$. Coordinates of the mapping $f$*

$$f(c) = c^l f^i_l \overline{e}_i \tag{9.2}$$

*relative to the basis $\overline{\overline{e}}$ satisfy to the equation*

$$f^j_k = f^i_0 C^j_{ki} \tag{9.3}$$

$$d = f^{i}_{0}\overline{e}_{i} \tag{9.4}$$

*Proof.* $(\mathrm{Re}\,D)$-linear mapping satisfies to equations

$$f(ac) = (ac)^{i} f^{j}_{i}\overline{e}_{j} = a^{k}c^{l}C^{i}_{kl}f^{j}_{i}\overline{e}_{j} \tag{9.5}$$

$$af(c) = a^{k}(f(c))^{i}C^{j}_{ki}\overline{e}_{j} = a^{k}c^{l}f^{i}_{l}C^{j}_{ki}\overline{e}_{j} \tag{9.6}$$

From equations (9.5), (9.6), it follows that

$$a^{k}c^{l}C^{i}_{kl}f^{j}_{i}\overline{e}_{j} = a^{k}c^{l}f^{i}_{l}C^{j}_{ki}\overline{e}_{j} \tag{9.7}$$

The equation

$$C^{i}_{kl}f^{j}_{i} = f^{i}_{l}C^{j}_{ki} \tag{9.8}$$

follows from the equation (9.7).

From the equation (9.8) it follows that[6]

$$C^{i}_{k0}f^{j}_{i} = f^{i}_{0}C^{j}_{ki} \tag{9.9}$$

From equations (8.13), (9.9), it follows that

$$\delta^{i}_{k}f^{j}_{i} = f^{i}_{0}C^{j}_{ki} \tag{9.10}$$

The equation (9.3) follows from the equation (9.10).

From equations (9.2), (9.3), it follows that

$$f(c) = c^{l}f^{j}_{0}C^{i}_{lj}\overline{e}_{i} \tag{9.11}$$

The equation (9.4) follows from the equations (9.1), (9.11). □

**Definition 9.2.** Let $D$ be the ring with conjugation. The mapping

$$f : D \rightarrow D$$

is called **$D\star$-antilinear**,[7] if the mapping $f$ satisfies to the equation

$$f(da) = f(a)d^{*} \tag{9.12}$$

□

**Theorem 9.3.** *$D\star$-antilinear mapping of the ring with conjugation $D$ has form*

$$f(c) = dc^{*} = dI \circ c \quad d \in D \tag{9.13}$$

*Let $\overline{\overline{e}}$ be $(\mathrm{Re}\,D)$-basis of the ring $D$. Coordinates of the mapping $f$*

$$f(c) = c^{l}f^{i}_{l}\overline{e}_{i} \tag{9.14}$$

*relative to the basis $\overline{\overline{e}}$ satisfy to the equation*

$$f^{j}_{l} = f^{k}_{0}I^{i}_{l}C^{j}_{ki} \tag{9.15}$$

$$d = f^{i}_{0}\overline{e}_{i} \tag{9.16}$$

---

[6] Let $l = 0$.

[7] I made the definition 9.2 by analogy with the definition 5.1. This definition also takes into account the requirements for antilinear homomorphism in the definition 10.2.

*Proof.* $(\mathrm{Re}\, D)$-linear mapping satisfies to equations

$$f(ac)=((ac)^*)^{i}f_{i}^{j}\overline{e}_{j}=c^{p}a^{q}(C_{pq}^{i})^*f_{i}^{j}\overline{e}_{j}\tag{9.17}$$

$$\begin{aligned}f(c)a^*&=(f(c))^{k}(a^*)^{i}C_{ki}^{j}\overline{e}_{j}=(c^*)^{l}f_{l}^{k}(a^*)^{i}C_{ki}^{j}\overline{e}_{j}\\&=I_{p}^{l}c^{p}f_{l}^{k}I_{q}^{i}a^{q}C_{ki}^{j}\overline{e}_{j}\end{aligned}\tag{9.18}$$

From equations (9.17), (9.18), it follows that

$$c^{p}a^{q}(C_{pq}^{i})^*f_{i}^{j}\overline{e}_{j}=I_{p}^{l}c^{p}f_{l}^{k}I_{q}^{i}a^{q}C_{ki}^{j}\overline{e}_{j}\tag{9.19}$$

The equation

$$(C_{pq}^{i})^*f_{i}^{j}=I_{p}^{l}f_{l}^{k}I_{q}^{i}C_{ki}^{j}\tag{9.20}$$

follows from the equation (9.19).

From the equation (9.20) it follows that[8]

$$(C_{0q}^{i})^*f_{i}^{j}=I_{0}^{l}f_{l}^{k}I_{q}^{i}C_{ki}^{j}\tag{9.21}$$

From equations (8.13), (8.19), (9.21), it follows that

$$\delta_{q}^{i}f_{i}^{j}=\delta_{0}^{l}f_{l}^{k}I_{q}^{i}C_{ki}^{j}\tag{9.22}$$

The equation (9.15) follows from the equation (9.22).

From equations (9.14), (9.15), it follows that

$$f(c)=c^{l}f_{0}^{k}I_{l}^{i}C_{ki}^{j}\overline{e}_{j}=f_{0}^{k}I_{l}^{i}c^{l}C_{ki}^{j}\overline{e}_{j}=f_{0}^{k}(c^*)^{i}C_{ki}^{j}\overline{e}_{j}\tag{9.23}$$

The equation (9.16) follows from the equations (9.1), (9.23). □

**Theorem 9.4.** *If mapping*

$$f:D\to D$$

*is $D\star$-linear mapping, then mapping*

$$g:D\to D\quad g(d)=(f(d))^*$$

*is $D\star$-antilinear mapping.*

*Proof.* The theorem follows from the equation

$$g(cd)=(f(cd))^*=(cf(d))^*=(f(d))^*\,c^*=g(d)\,c^*$$

and the definition 9.2. □

## 10. Antilinear Homomorphism of $D\star$-Algebra

In this section $D$ is a ring with conjugation.

**Theorem 10.1.** *Let diagram of representations*

(10.1)

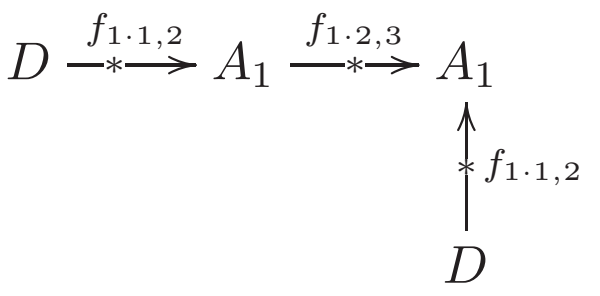

$$f_{1\cdot 1,2}(d):\overline{v}\to d\,\overline{v}$$
$$f_{1\cdot 2,3}(\overline{v}):\overline{w}\to \overline{C}_1(\overline{v},\overline{w})$$
$$\overline{C}_1\in\mathcal{L}(A_1^2;A_1)$$

[8] Let $p=0$.

*describe $D\star$-algebra $A_1$. Let diagram of representations*

$$D \xrightarrow[*]{f_{2\cdot 1,2}} A_2 \xrightarrow[*]{f_{2\cdot 2,3}} A_2 \qquad \uparrow^{*} f_{2\cdot 1,2} \quad D \tag{10.2}$$

$$f_{2\cdot 1,2}(d) : \overline{v} \to \overline{v}\,d$$
$$f_{2\cdot 2,3}(\overline{v}) : \overline{w} \to \overline{C}_2(\overline{w},\overline{v})$$
$$\overline{C}_2 \in \mathcal{L}(A_2^2; A_2)$$

*describe $\star D$-algebra $A_2$. Morphism of $D\star$-algebra $A_1$ into $\star D$-algebra $A_2$*

$$r_1 : D \to D \quad \overline{r}_2 : A_1 \to A_2$$
$$r_1(d) = d^*$$

*is $D\star$-antilinear mapping of $D\star$-algebra $A_1$ into $\star D$-algebra $A_2$ such that*

$$\overline{r}_2(\overline{a}\overline{b}) = \overline{r}_2(\overline{b})\overline{r}_2(\overline{a}) \tag{10.3}$$

*Proof.* According to equation [5]-(4.2.3), morphism $(r_1, \overline{r}_2)$ of representation $f_{1,2}$ satisfies to the equation

$$\begin{aligned}\overline{r}_2(f_{1\cdot 1,2}(d)(\overline{a})) &= f_{2\cdot 1,2}(r_1(d))(\overline{r}_2(\overline{a}))\\ \overline{r}_2(d\,\overline{a}) &= \overline{r}_2(\overline{a})\,d^*\end{aligned} \tag{10.4}$$

Therefore, the mapping $(r_1, \overline{r}_2)$ is $D\star$-antilinear mapping.

According to equation [5]-(4.2.3), the morphism $(\overline{r}_2, \overline{r}_2)$ of representation $f_{2,3}$ satisfies to the equation[9]

$$\overline{r}_2(f_{1\cdot 2,3}(\overline{a}_2)(\overline{a}_3)) = f_{2\cdot 2,3}(\overline{r}_2(\overline{a}_2))(\overline{r}_2(\overline{a}_3)) \tag{10.5}$$

From equations (10.5), (10.1), (10.2), it follows that

$$\overline{r}_2(\overline{C}_1(\overline{v},\overline{w})) = \overline{C}_2(\overline{r}_2(\overline{w}), \overline{r}_2(\overline{v})) \tag{10.6}$$

Equation (10.3) follows from equations (10.6), (6.3). □

**Definition 10.2.** The morphism of representations of $D\star$-algebra $A_1$ into $\star D$-algebra $A_2$ is called **$D\star$-antilinear homomorphism** of $D\star$-algebra $A_1$ into $\star D$-algebra $A_2$. □

**Theorem 10.3.** *Let $\overline{\overline{e}}_1$ be the basis of $D\star$-algebra $A_1$. Let $\overline{\overline{e}}_2$ be the basis of $\star D$-algebra $A_2$. Then antilinear homomorphism[10] $(r_1, \overline{r}_2)$ of $D\star$-algebra $A_1$ into $\star D$-algebra $A_2$ has presentation[11]*

$$\begin{aligned}\overline{b} &= \overline{e}_{2*}{}^*r_{2*}{}^*r_1(a) = \overline{e}_{2\cdot i}\, r_{2\cdot}{}^{i}_{j}\, r_1(a^j)\\ b &= r_{2*}{}^*r_1(a)\end{aligned} \tag{10.7}$$

*relative to selected bases. Here*

- *$a$ is coordinate matrix of vector $\overline{a}$ relative the basis $\overline{\overline{e}}_1$.*
- *$b$ is coordinate matrix of vector*

$$\overline{b} = \overline{r}_2(\overline{a})$$

*relative the basis $\overline{\overline{e}}_2$.*

---

[9]Since in diagrams of representations (6.4), (6.5), supports of $\Omega_2$-algebra and $\Omega_3$-algebra coincide, then morphisms of representations on levels 2 and 3. coincide also

[10]This theorem is similar to the theorem [2]-4.4.3.

[11]For complex field, the equation (5.11) is consequence of the equation (10.7).

- *$r_2$ is coordinate matrix of set of vectors $(\overline{r}_2(\overline{e}_{1\cdot i}))$ relative the basis $\overline{\overline{e}}_2$. The matrix $r_2$ is called* **matrix of antilinear homomorphism** *relative bases $\overline{\overline{e}}_1$ and $\overline{\overline{e}}_2$.*

*Proof.* Vector $\overline{a} \in A_1$ has expansion

$$\overline{a} = a^*{}_*\overline{e}_1$$

relative to the basis $\overline{\overline{e}}_1$. Vector $\overline{b} \in A_2$ has expansion

$$\overline{b} = \overline{e}_{2*}{}^*b \tag{10.8}$$

relative to the basis $\overline{\overline{e}}_2$.

Since $(r_1, \overline{r}_2)$ is an antilinear homomorphism, then from (10.4), it follows that

$$\overline{b} = \overline{r}_2(\overline{a}) = \overline{r}_2(a^*{}_*\overline{e}_1) = \overline{r}_2(\overline{e}_1)_*{}^*a^* \tag{10.9}$$

where

$$a^* = \begin{pmatrix} (a^1)^* \\ ... \\ (a^n)^* \end{pmatrix}$$

$\overline{r}_2(\overline{e}_{1\cdot i})$ is also a vector of $\star D$-module $A_2$ and has expansion

$$\overline{r}_2(\overline{e}_{1\cdot i}) = \overline{e}_{2*}{}^*r_{2\cdot i} = \overline{e}_j\, r_{2\cdot}{}^j_i \tag{10.10}$$

relative to basis $\overline{\overline{e}}_2$. Combining (10.9) and (10.10) we get

$$\overline{b} = \overline{e}_{2\circ}{}^\circ r_{2\circ}{}^\circ a^* \tag{10.11}$$

(10.7) follows from comparison of (10.8) and (10.11) and theorem [2]-4.3.3. □

**Theorem 10.4.** *Let $\overline{\overline{e}}_1$ be the basis of $D\star$-algebra $A_1$. Let $\overline{\overline{e}}_2$ be the basis of $\star D$-algebra $A_2$. There exists relation between the matrix of antilinear homomorphism and structural constants*

$$(C_{1\cdot}{}^k_{ij})^*\, r_{2\cdot}{}^l_k = r_{2\cdot}{}^p_j\, r_{2\cdot}{}^q_i\, C_{2\cdot}{}^l_{pq} \tag{10.12}$$

*Proof.* Let

$$\overline{a}, \overline{b} \in A_1 \quad \overline{a} = a^*{}_*\overline{e}_1 \quad \overline{b} = b^*{}_*\overline{e}_1$$

From equations (6.2), (6.3), (10.1), it follows that

$$\overline{a}\overline{b} = a^i b^j C_{1\cdot}{}^k_{ij}\overline{e}_{1\cdot k} \tag{10.13}$$

From equations (10.4), (10.13), it follows that

$$\overline{r}_2(\overline{a}\overline{b}) = \overline{r}_2(\overline{e}_{1\cdot k})(a^i b^j C_{1\cdot}{}^k_{ij})^* = \overline{r}_2(\overline{e}_{1\cdot k})(C_{1\cdot}{}^k_{ij})^*(b^j)^*(a^i)^* \tag{10.14}$$

From the theorem 10.3 and the equation (10.14), it follows that

$$\overline{r}_2(\overline{a}\overline{b}) = \overline{e}_{2\cdot l}\, r_{2\cdot}{}^l_k(C_{1\cdot}{}^k_{ij})^*(a^i)^*(b^j)^* \tag{10.15}$$

From the equation (10.3) and the theorem 10.3, it follows that

$$\overline{r}_2(\overline{a}\overline{b}) = \overline{r}_2(\overline{b})\overline{r}_2(\overline{a}) = \overline{e}_{2\cdot q}\, r_{2\cdot}{}^q_j(b^j)^*\, \overline{e}_{2\cdot p}\, r_{2\cdot}{}^p_i(a^i)^* \tag{10.16}$$

From equations (6.2), (6.3), (10.2), (10.16), it follows that

$$\overline{r}_2(\overline{a}\overline{b}) = \overline{e}_{2\cdot l}\, C_{2\cdot}{}^l_{pq} r_{2\cdot}{}^p_i(a^i)^* r_{2\cdot}{}^q_j(b^j)^* \tag{10.17}$$

From equations (10.15), (10.17), it follows that

$$\overline{e}_{2\cdot l}\, r_{2\cdot}{}^{l}_{k}(C_{1\cdot}{}^{k}_{ij})^*(a^i)^*(b^j)^* = \overline{e}_{2\cdot l}\, C_{2\cdot}{}^{l}_{pq} r_{2\cdot}{}^{p}_{i}(a^i)^* r_{2\cdot}{}^{q}_{j}(b^j)^* \tag{10.18}$$

The equation (10.12) follows from the equation (10.18), because vectors of the basis $\overline{\overline{e}}_2$ are linear indepensent, and $a^i$, $b^i$ (and therefore, $r_1(a^i)$, $r_1(b^i)$) are arbitrary values. □

## 11. Involution of Banach $D$-algebra

In this section $D$ is a ring with conjugation. We also will assume that there is structure of $D\star$-module and $\star D$-module in Abelian group $A$.[12]

**Definition 11.1.** Let $A$ be a Banach $D$-algebra. Involution is an antilinear homomorphism $x \rightarrow x^*$, $x \in A$, satisfying

$$x^{**} = x \tag{11.1}$$

$$\|x^*\,x\| = \|x\|^2$$

$$x, y \in A$$

□

**Theorem 11.2.** *Let $A$ be a Banach $D$-algebra. Let mapping $x \rightarrow x^*$, $x \in A$, be involution. Then*

$$(xy)^* = y^* x^* \tag{11.2}$$

*Proof.* The equation (11.2) is consequence of the definition 11.1 and the theorem 10.1. □

**Definition 11.3.** The basis $\overline{\overline{e}}$ is called normal basis, if $\|\overline{e}_i\| = 1$ for any vector $\overline{e}_i$ of the basis $\overline{\overline{e}}$. □

Without loss of generality, we can assume that the basis $\overline{\overline{e}}$ is normal basis. If we assume that the norm of the vector $\overline{e}_i$ is different from 1, then we can replace this vector by the vector

$$\overline{e}'_i = \frac{1}{\|\overline{e}_i\|}\overline{e}_i$$

**Theorem 11.4.** *Let $D$ be commutative ring. If the mapping*

$$f : A_1 \rightarrow A_2$$

*is linear mapping, then the mapping*

$$h(x) = f(x^*)$$

*is antilinear mapping. If the mapping*

$$f : A_1 \rightarrow A_2$$

*is linear homomorphism, then the mapping*

$$h(x) = f(x^*)$$

*is antilinear homomorphism.*

[12] Conditions when this requirement is fulfilled, will be discussed in a separate paper. This condition is evident if $D$ is field. This condition is satisfied according to the theorem [2]-6.3.1 if $D$ is division ring.

*Proof.* According to the theorem 10.1 and the equation (10.4)

$$(da)^* = a^*\,d^* \quad d \in D \quad a \in A \tag{11.3}$$

If $f$ is linear mapping, then

$$f(dx) = df(x) \tag{11.4}$$

From equations (11.3), (11.4), it follows that

$$h(xd) = f((xd)^*) = f(d^*\,x^*) = d^*\,f(x^*) = d^*\,h(x)$$

If $f$ is linear homomorphism, then from equations (6.6), (11.2), it follows that

$$h(ab) = f((ab)^*) = f(b^*\,a^*) = f(b^*)f(a^*) = h(b)h(a)$$

Therefore, $h$ is antilinear homomorphism. □

## 12. $C^*$-algebra

**Definition 12.1.** A $C^*$**-algebra** $A$ is a Banach $C$-algebra that is endowed with an involution.[13] □

According to the theorem 5.4, involution has form

$$\begin{aligned}(x^*)^{0j} &= f^{0j}_{0i}x^{0i} + f^{1j}_{0i}x^{1i}\\ (x^*)^{1j} &= f^{1j}_{0i}x^{0i} - f^{0j}_{0i}x^{1i}\end{aligned} \tag{12.1}$$

relative to the basis $\overline{\overline{e}}_{AR}$ .

**Theorem 12.2.**

$$\begin{aligned}\delta^j_k &= f^{0j}_{0i}f^{0i}_{0k} + f^{1j}_{0i}f^{1i}_{0k}\\ 0 &= f^{0j}_{0i}f^{1i}_{0k} - f^{1j}_{0i}f^{0i}_{0k}\end{aligned} \tag{12.2}$$

*Proof.* From the equation (12.1) it follows that

$$\begin{aligned}(x^{**})^{0j} &= f^{0j}_{0i}(x^*)^{0i} + f^{1j}_{0i}(x^*)^{1i}\\ &= f^{0j}_{0i}(f^{0i}_{0k}x^{0k} + f^{1i}_{0k}x^{1k}) + f^{1j}_{0i}(f^{1i}_{0k}x^{0k} - f^{0i}_{0k}x^{1k})\\ &= f^{0j}_{0i}f^{0i}_{0k}x^{0k} + f^{0j}_{0i}f^{1i}_{0k}x^{1k} + f^{1j}_{0i}f^{1i}_{0k}x^{0k} - f^{1j}_{0i}f^{0i}_{0k}x^{1k}\\ (x^{**})^{1j} &= f^{1j}_{0i}(x^*)^{0i} - f^{0j}_{0i}(x^*)^{1i}\\ &= f^{1j}_{0i}(f^{0i}_{0k}x^{0k} + f^{1i}_{0k}x^{1k}) - f^{0j}_{0i}(f^{1i}_{0k}x^{0k} - f^{0i}_{0k}x^{1k})\\ &= f^{1j}_{0i}f^{0i}_{0k}x^{0k} + f^{1j}_{0i}f^{1i}_{0k}x^{1k} - f^{0j}_{0i}f^{1i}_{0k}x^{0k} + f^{0j}_{0i}f^{0i}_{0k}x^{1k}\end{aligned} \tag{12.3}$$

(12.2) follows from equations (11.1), (12.3). □

**Theorem 12.3.**

$$\begin{aligned}&C_{AC\cdot}{}^{k0}_{ij}f^{0l}_{0k} + C_{AC\cdot}{}^{k1}_{ij}f^{1l}_{0k}\\ &= f^{0p}_{0j}f^{0q}_{0i}C_{AC\cdot}{}^{l0}_{pq} - f^{1p}_{0j}f^{1q}_{0i}C_{AC\cdot}{}^{l0}_{pq} - f^{0p}_{0j}f^{1q}_{0i}C_{AC\cdot}{}^{l1}_{pq} - f^{1p}_{0j}f^{0q}_{0i}C_{AC\cdot}{}^{l1}_{pq}\\ &C_{AC\cdot}{}^{k0}_{ij}f^{1l}_{0k} - C_{AC\cdot}{}^{k1}_{ij}f^{0l}_{0k}\\ &= f^{0p}_{0j}f^{0q}_{0i}C_{AC\cdot}{}^{l1}_{pq} - f^{1p}_{0j}f^{1q}_{0i}C_{AC\cdot}{}^{l1}_{pq} + f^{0p}_{0j}f^{1q}_{0i}C_{AC\cdot}{}^{l0}_{pq} + f^{1p}_{0j}f^{0q}_{0i}C_{AC\cdot}{}^{l0}_{pq}\end{aligned} \tag{12.4}$$

[13] This definition is based on the definition [6]-2.2.1.

*Proof.* From the equation (5.11) it follows that

$$f_i^j = f_{0i}^{0j} + f_{0i}^{1j}i \tag{12.5}$$

From equations (3.17), (12.5), (10.12), it follows that

$$\begin{aligned}&(C_{AC\cdot}{}_{ij}^{k})^* f_k^l = f_j^p f_i^q C_{AC\cdot}{}_{pq}^{l}\\&(C_{AC\cdot}{}_{ij}^{k0} + C_{AC\cdot}{}_{ij}^{k1}i)^*(f_{0k}^{0l} + f_{0k}^{1l}i)\\=&(f_{0j}^{0p} + f_{0j}^{1p}i)(f_{0i}^{0q} + f_{0i}^{1q}i)(C_{AC\cdot}{}_{pq}^{l0} + C_{AC\cdot}{}_{pq}^{l1}i)\end{aligned} \tag{12.6}$$

From the equation (12.6) it follows that

$$\begin{aligned}&C_{AC\cdot}{}_{ij}^{k0} f_{0k}^{0l} + C_{AC\cdot}{}_{ij}^{k1} f_{0k}^{1l} + (C_{AC\cdot}{}_{ij}^{k0} f_{0k}^{1l} - C_{AC\cdot}{}_{ij}^{k1} f_{0k}^{0l})i\\=&(f_{0j}^{0p} f_{0i}^{0q} - f_{0j}^{1p} f_{0i}^{1q} + (f_{0j}^{0p} f_{0i}^{1q} + f_{0j}^{1p} f_{0i}^{0q})i)(C_{AC\cdot}{}_{pq}^{l0} + C_{AC\cdot}{}_{pq}^{l1}i)\\=&f_{0j}^{0p} f_{0i}^{0q} C_{AC\cdot}{}_{pq}^{l0} - f_{0j}^{1p} f_{0i}^{1q} C_{AC\cdot}{}_{pq}^{l0} - f_{0j}^{0p} f_{0i}^{1q} C_{AC\cdot}{}_{pq}^{l1} - f_{0j}^{1p} f_{0i}^{0q} C_{AC\cdot}{}_{pq}^{l1}\\&+(f_{0j}^{0p} f_{0i}^{0q} C_{AC\cdot}{}_{pq}^{l1} - f_{0j}^{1p} f_{0i}^{1q} C_{AC\cdot}{}_{pq}^{l1} + f_{0j}^{0p} f_{0i}^{1q} C_{AC\cdot}{}_{pq}^{l0} + f_{0j}^{1p} f_{0i}^{0q} C_{AC\cdot}{}_{pq}^{l0})i\end{aligned} \tag{12.7}$$

The equation (12.4) follows from the equation (12.7). □

## 13. Example of Involution

Let

$$e_1^1 = \begin{pmatrix}1 & 0\\0 & 0\end{pmatrix} \quad e_1^2 = \begin{pmatrix}0 & 1\\0 & 0\end{pmatrix}$$
$$e_2^1 = \begin{pmatrix}0 & 0\\1 & 0\end{pmatrix} \quad e_2^2 = \begin{pmatrix}0 & 0\\0 & 1\end{pmatrix}$$

be basis $\overline{\overline{e}}$ of algebra of $2 \times 2$ matrices. Matrix

$$a = \begin{pmatrix}a_1^1 & a_2^1\\a_1^2 & a_2^2\end{pmatrix}$$

has expansion

$$a = a_1^1 e_1^1 + a_2^1 e_1^2 + a_1^2 e_2^1 + a_2^2 e_2^2$$

relative to the basis $\overline{\overline{e}}$. Structural constants in algebra of matrices have form

$$C^{\cdot p}_{q}{}^{\cdot}{}^{i}_{j}{}^{\cdot}{}^{k}_{l} = \delta_l^i \delta_j^p \delta_q^k$$

The product of matrices $a$ and $b$ has form

$$(ab)_q^p = C^{\cdot p}_{q}{}^{\cdot}{}^{i}_{j}{}^{\cdot}{}^{k}_{l} a_i^j b_k^l = \delta_l^i \delta_j^p \delta_q^k a_i^j b_k^l = a_l^p b_q^l$$

If we represent the matrix $a$ as vector

$$a = \begin{pmatrix}a_1^1\\a_2^1\\a_1^2\\a_2^2\end{pmatrix}$$

then we can represent the matrix of transformation $f$

$$a'^i_j = f^{\cdot i}_{j\cdot l}{}^k a^l_k$$

in the following form

$$f = \begin{pmatrix} f^{\cdot 1}_{1\cdot 1}{}^1 & f^{\cdot 1}_{1\cdot 2}{}^1 & f^{\cdot 1}_{1\cdot 1}{}^2 & f^{\cdot 1}_{1\cdot 2}{}^2 \\ f^{\cdot 1}_{2\cdot 1}{}^1 & f^{\cdot 1}_{2\cdot 2}{}^1 & f^{\cdot 1}_{2\cdot 1}{}^2 & f^{\cdot 1}_{2\cdot 2}{}^2 \\ f^{\cdot 2}_{1\cdot 1}{}^1 & f^{\cdot 2}_{1\cdot 2}{}^1 & f^{\cdot 2}_{1\cdot 1}{}^2 & f^{\cdot 2}_{1\cdot 2}{}^2 \\ f^{\cdot 2}_{2\cdot 1}{}^1 & f^{\cdot 2}_{2\cdot 2}{}^1 & f^{\cdot 2}_{2\cdot 1}{}^2 & f^{\cdot 2}_{2\cdot 2}{}^2 \end{pmatrix}$$

The matrix

$$f = \begin{pmatrix} 1 & 0 & 0 & 0 \\ 0 & 0 & 1 & 0 \\ 0 & 1 & 0 & 0 \\ 0 & 0 & 0 & 1 \end{pmatrix}$$

corresponds to the transposition

$$a = \begin{pmatrix} a^1_1 & a^1_2 \\ a^2_1 & a^2_2 \end{pmatrix} \to a^T = \begin{pmatrix} a^1_1 & a^2_1 \\ a^1_2 & a^2_2 \end{pmatrix}$$

Since we do not have analytic representation for transposition operator, then we use nonstandard format of representation.

$$\begin{aligned}
(a^T b^T)^p_q &= C^{\cdot p}_{q\cdot j}{}^i{}_{\cdot l}{}^k (a^T)^j_i (b^T)^l_k \\
&= \delta^i_l \delta^p_j \delta^k_q \begin{pmatrix} 1 & 0 & 0 & 0 \\ 0 & 0 & 1 & 0 \\ 0 & 1 & 0 & 0 \\ 0 & 0 & 0 & 1 \end{pmatrix} \begin{pmatrix} a^1_1 \\ a^1_2 \\ a^2_1 \\ a^2_2 \end{pmatrix} \begin{pmatrix} 1 & 0 & 0 & 0 \\ 0 & 0 & 1 & 0 \\ 0 & 1 & 0 & 0 \\ 0 & 0 & 0 & 1 \end{pmatrix} \begin{pmatrix} b^1_1 \\ b^1_2 \\ b^2_1 \\ b^2_2 \end{pmatrix} \\
&= \delta^i_l \delta^p_j \delta^k_q \begin{pmatrix} {}^j_i : {}^1_1 & a^1_1 \\ {}^j_i : {}^1_2 & a^2_1 \\ {}^j_i : {}^2_1 & a^1_2 \\ {}^j_i : {}^2_2 & a^2_2 \end{pmatrix} \begin{pmatrix} {}^l_k : {}^1_1 & b^1_1 \\ {}^l_k : {}^1_2 & b^2_1 \\ {}^l_k : {}^2_1 & b^1_2 \\ {}^l_k : {}^2_2 & b^2_2 \end{pmatrix} = \begin{pmatrix} {}^p_q : {}^1_1 & a^1_1 b^1_1 + a^2_1 b^1_2 \\ {}^p_q : {}^1_2 & a^1_1 b^2_1 + a^2_1 b^2_2 \\ {}^p_q : {}^2_1 & a^1_2 b^1_1 + a^2_2 b^1_2 \\ {}^p_q : {}^2_2 & a^1_2 b^2_1 + a^2_2 b^2_2 \end{pmatrix} \qquad (13.1) \\
&= \begin{pmatrix} 1 & 0 & 0 & 0 \\ 0 & 0 & 1 & 0 \\ 0 & 1 & 0 & 0 \\ 0 & 0 & 0 & 1 \end{pmatrix} \begin{pmatrix} b^1_k a^k_1 \\ b^1_k a^k_2 \\ b^2_k a^k_1 \\ b^2_k a^k_2 \end{pmatrix} \\
&= ((ba)^T)^p_q
\end{aligned}$$

If $a^i_j$ are complex numbers, then involution can be represented by antilinear mapping with matrix

$$f=\begin{pmatrix}I&0&0&0\\0&0&I&0\\0&I&0&0\\0&0&0&I\end{pmatrix}=\left(\begin{array}{cc|cc|cc|cc}1&0&0&0&0&0&0&0\\0&-1&0&0&0&0&0&0\\\hline 0&0&0&0&1&0&0&0\\0&0&0&0&0&-1&0&0\\\hline 0&0&1&0&0&0&0&0\\0&0&0&-1&0&0&0&0\\\hline 0&0&0&0&0&0&1&0\\0&0&0&0&0&0&0&-1\end{array}\right)$$

## 14. References

[1] Serge Lang, Algebra, Springer, 2002
[2] Aleks Kleyn, Lectures on Linear Algebra over Division Ring, eprint arXiv:math.GM/0701238 (2010)
[3] Aleks Kleyn, Quaternion Rhapsody, eprint arXiv:0909.0855 (2010)
[4] Aleks Kleyn, Linear Mappings of Free Algebra: First Steps in Noncommutative Linear Algebra, Lambert Academic Publishing, 2010
[5] Aleks Kleyn, Representation Theory: Representation of Universal Algebra, Lambert Academic Publishing, 2011
[6] William Arveson, A short course on spectral theory. Springer - Verlag, New York, 2002
[7] Robert Hermann, Topics in the mathematics of quantum mechanics. Math Sci Press, 1973

## 15. Index

## 16. Special Symbols and Notations

# $C^*$-рапсодия

Александр Клейн

Аннотация. В статье рассмотрены линейные и антилинейные отображения конечномерной алгебры над комплексным полем, а также отображение инволюции. Рассмотрен также пример $C^*$-алгебры.

## Содержание

## 1. Предисловие

Поскольку я изучаю математический анализ в банаховой алгебре, трудно обойти вниманием $C^*$-алгебру, которая является инструментом в квантовой механике. Однако, прежде чем рассмотреть конкретные результаты, мне нужно было понять структуру линейного отображения $C^*$-алгебры. Это послужило причиной написать эту статью.

Однако существуют различные алгебры, где определено сопряжение. Поэтому я решил не ограничивать себя полем комплексных чисел. Рассмотрение

Aleks_Kleyn@MailAPS.org.
http://sites.google.com/site/AleksKleyn/.
http://arxiv.org/a/kleyn_a_1.
http://AleksKleyn.blogspot.com/.

линейных и антилинейных отображений алгебры с сопряжением оказалось отличным от концепции, которую я рассматривал в [4]. Заманчиво найти утверждение, аналогично тому, что аддитивное отображение поля комплексных чисел можно разложить в сумму линейного и антилинейного отображений. Подобное утверждение оказалось бы мостиком между результатами в [4] и этой статьёй. Однако есть основания полагать, что в общем случае это утверждение неверно.

Это означает, что в общем случае рассмотрение сопряжения недостаточно для полной классификации линейных отображений. Однако эта классификация важна для того, чтобы понять структуру замкнутой дифференциальной формы в банаховой алгебре.

## 2. Соглашения

(1) Пусть $A$ - свободная конечно мерная алгебра. При разложении элемента алгебры $A$ относительно базиса $\overline{\overline{e}}$ мы пользуемся одной и той же корневой буквой для обозначения этого элемента и его координат. Однако в алгебре не принято использовать векторные обозначения. В выражении $a^2$ не ясно - это компонента разложения элемента $a$ относительно базиса или это операция возведения в степень. Для облегчения чтения текста мы будем индекс элемента алгебры выделять цветом. Например,

$$a = a^{i}\overline{e}_{i}$$

(2) Если свободная конечномерная алгебра имеет единицу, то мы будем отождествлять вектор базиса $\overline{e}_{0}$ с единицей алгебры.

(3) Хотя алгебра является свободным модулем над некоторым кольцом, мы не пользуемся векторными обозначениями при записи элементов алгебры. В тех случаях, когда я рассматриваю матрицу координат элемента алгебры, я буду пользоваться векторными обозначениями для записи соответствующего элемента. Чтобы не возникала неоднозначность при записи сопряжения, я буду обозначать $a^*$ элемент, сопряжённый элементу $a$.

(4) Без сомнения, у читателя могут быть вопросы, замечания, возражения. Я буду признателен любому отзыву.

## 3. $C$-алгебра

Рассматриваемое ниже построение основано на построении, выполненном в разделе [3]-5. В этой статье мы будем пользоваться обозначениями, принятыми в этом разделе.

Пусть $C$ - поле комплексных чисел. Пусть $\overline{\overline{e}}_C$

$$\overline{e}_{C\cdot 0} = 1 \quad \overline{e}_{C\cdot 1} = i \tag{3.1}$$

базис алгебры $C$ над полем действительных чисел $R$. Произведение в алгебре $C$ определено согласно правилу

$$\overline{e}_{C\cdot 1}^{2} = -\overline{e}_{C\cdot 0} \tag{3.2}$$

Согласно теореме [3]-3.1 структурные константы поля комплексных чисел $C$ над полем действительных чисел $R$ имеют вид

$$\begin{array}{ll} C_{C\cdot}{}^{0}_{00}=1 & C_{C\cdot}{}^{1}_{01}= \ \ 1 \\ C_{C\cdot}{}^{1}_{10}=1 & C_{C\cdot}{}^{0}_{11}=-1 \end{array} \tag{3.3}$$

Согласно определению [4]-2.2.1, $C$-алгебра $A$ является $C$-векторным пространством. Пусть $\overline{\overline{e}}_{AC}$ - $C$-базис конечно мерной $C$-алгебры $A$. Произведение в $C$-алгебре $A$ определено равенством

$$\overline{e}_{AC\cdot k}\overline{e}_{AC\cdot l}=C_{AC\cdot}{}^{i}_{kl}\overline{e}_{AC\cdot i} \tag{3.4}$$

где $C_{AC\cdot}{}^{k}_{ij}$ - структурные константы $C$-алгебры $A$ над полем $C$.

Я буду рассматривать $C$-алгебру $A$ как прямую сумму алгебр $C$. Каждое слагаемое суммы я отождествляю с вектором базиса $\overline{\overline{e}}_{AC}$ . Соответственно, я могу рассматривать $C$-алгебру $A$ как алгебру над полем $R$. Пусть $\overline{\overline{e}}_{AR}$ - базис $C$-алгебры $A$ над полем $R$. Индекс базиса $\overline{\overline{e}}_{AR}$ состоит из двух индексов: индекса слоя и индекса вектора базиса $\overline{\overline{e}}_{C}$ в слое.

Я буду отождествлять вектор базиса $\overline{e}_{AC\cdot i}$ с единицей в соответствующем слое. Тогда

$$\overline{e}_{AR\cdot ji}=\overline{e}_{C\cdot j}\overline{e}_{AC\cdot i} \tag{3.5}$$

Произведение векторов базиса $\overline{\overline{e}}_{AR}$ имеет вид

$$\overline{e}_{AR\cdot ji}\overline{e}_{AR\cdot mk}=\overline{e}_{C\cdot j}\overline{e}_{AC\cdot i}\overline{e}_{C\cdot m}\overline{e}_{AC\cdot k}=C_{C\cdot}{}^{a}_{jm}\overline{e}_{C\cdot a}C_{AC\cdot}{}^{b}_{ik}\overline{e}_{AC\cdot b} \tag{3.6}$$

$$\overline{e}_{AR\cdot 00}\,\overline{e}_{AR\cdot ki}=\overline{e}_{AR\cdot ki} \tag{3.7}$$

$$\overline{e}_{AR\cdot 10}\,\overline{e}_{AR\cdot 1i}=-\overline{e}_{AR\cdot 0i} \tag{3.8}$$

Так как $C_{AC\cdot}{}^{b}_{ik}\in C$ , то разложение $C_{AC\cdot}{}^{b}_{ik}\in C$ относительно базиса $\overline{\overline{e}}_{C}$ имеет вид

$$C_{AC\cdot}{}^{b}_{ik}=C_{AC\cdot}{}^{bc}_{ik}\overline{e}_{C\cdot c} \tag{3.9}$$

Мы можем представить преобразование $f$ с помощью матрицы $f^{ij}_{kl}$

$$a'^{ij}=f^{ij}_{kl}a^{kl} \tag{3.10}$$

$$\begin{array}{l} a'^{0j}=f^{0j}_{0l}a^{0l}+f^{0j}_{1l}a^{1l} \\ a'^{1j}=f^{1j}_{0l}a^{0l}+f^{1j}_{1l}a^{1l} \end{array} \tag{3.11}$$

Мы можем представить $c\in C$ в следующем виде

$$c=c^{0}\,\overline{e}_{AR\cdot 00}+c^{1}\,\overline{e}_{AR\cdot 10} \tag{3.12}$$

$$c^{*}=c^{0}\,\overline{e}_{AR\cdot 00}-c^{1}\,\overline{e}_{AR\cdot 10} \tag{3.13}$$

Из равенств (3.7), (3.8), (3.12), следует

$$(ca)^{0i}=c^{0}a^{0i}-c^{1}a^{1i} \tag{3.14}$$

$$(ca)^{1i}=c^{0}a^{1i}+c^{1}a^{0i} \tag{3.15}$$

**Теорема 3.1.** *Структурные константы $C$-алгебры $A$ над полем действительных чисел $R$ имеют вид*

$$
\left\{
\begin{aligned}
&C_{AR\cdot\cdot 0i\cdot 0k}^{\cdot 0b} = C_{AC\cdot ik}^{b0}\\
&C_{AR\cdot\cdot 1i\cdot 1k}^{\cdot 0b} = -C_{AC\cdot ik}^{b0}\\
&C_{AR\cdot\cdot 0i\cdot 1k}^{\cdot 0b} = -C_{AC\cdot ik}^{b1}\\
&C_{AR\cdot\cdot 1i\cdot 0k}^{\cdot 0b} = -C_{AC\cdot ik}^{b1}\\
&C_{AR\cdot\cdot 0i\cdot 1k}^{\cdot 1b} = C_{AC\cdot ik}^{b0}\\
&C_{AR\cdot\cdot 1i\cdot 0k}^{\cdot 1b} = C_{AC\cdot ik}^{b0}\\
&C_{AR\cdot\cdot 0i\cdot 0k}^{\cdot 1b} = C_{AC\cdot ik}^{b1}\\
&C_{AR\cdot\cdot 1i\cdot 1k}^{\cdot 1b} = -C_{AC\cdot ik}^{b1}
\end{aligned}
\right. \tag{3.16}
$$

*Доказательство.* Согласно равенствам [3]-(5.3), [3]-(5.5), структурные константы алгебры $A$ над полем $R$ имеют вид

$$
\begin{gathered}
C_{AR\cdot\cdot ji\cdot mk}^{\cdot db} = C_{C\cdot jm}^{a} C_{C\cdot ac}^{d} C_{AC\cdot ik}^{bc}\\
C_{AC\cdot ik}^{b} = C_{AC\cdot ik}^{bc} \overline{e}_{C\cdot c} = C_{AC\cdot ik}^{b0} + C_{AC\cdot ik}^{b1} i
\end{gathered} \tag{3.17}
$$

Из равенств (3.3), (3.17) следует

$$
\left\{
\begin{aligned}
C_{AR\cdot\cdot ji\cdot mk}^{\cdot 0b} &= C_{C\cdot jm}^{0} C_{C\cdot 00}^{0} C_{AC\cdot ik}^{b0} + C_{C\cdot jm}^{1} C_{C\cdot 11}^{0} C_{AC\cdot ik}^{b1}\\
&= C_{C\cdot jm}^{0} C_{AC\cdot ik}^{b0} - C_{C\cdot jm}^{1} C_{AC\cdot ik}^{b1}\\
C_{AR\cdot\cdot ji\cdot mk}^{\cdot 1b} &= C_{C\cdot jm}^{1} C_{C\cdot 10}^{1} C_{AC\cdot ik}^{b0} + C_{C\cdot jm}^{0} C_{C\cdot 01}^{1} C_{AC\cdot ik}^{b1}\\
&= C_{C\cdot jm}^{1} C_{AC\cdot ik}^{b0} + C_{C\cdot jm}^{0} C_{AC\cdot ik}^{b1}
\end{aligned}
\right. \tag{3.18}
$$

Из равенств (3.3), (3.18) следует

$$
\left\{
\begin{aligned}
&C_{AR\cdot\cdot 0i\cdot 0k}^{\cdot 0b} = C_{C\cdot 00}^{0} C_{AC\cdot ik}^{b0} - C_{C\cdot 00}^{1} C_{AC\cdot ik}^{b1}\\
&C_{AR\cdot\cdot 1i\cdot 1k}^{\cdot 0b} = C_{C\cdot 11}^{0} C_{AC\cdot ik}^{b0} - C_{C\cdot 11}^{1} C_{AC\cdot ik}^{b1}\\
&C_{AR\cdot\cdot 0i\cdot 1k}^{\cdot 0b} = C_{C\cdot 01}^{0} C_{AC\cdot ik}^{b0} - C_{C\cdot 01}^{1} C_{AC\cdot ik}^{b1}\\
&C_{AR\cdot\cdot 1i\cdot 0k}^{\cdot 0b} = C_{C\cdot 10}^{0} C_{AC\cdot ik}^{b0} - C_{C\cdot 10}^{1} C_{AC\cdot ik}^{b1}\\
&C_{AR\cdot\cdot 0i\cdot 1k}^{\cdot 1b} = C_{C\cdot 01}^{1} C_{AC\cdot ik}^{b0} + C_{C\cdot 01}^{0} C_{AC\cdot ik}^{b1}\\
&C_{AR\cdot\cdot 1i\cdot 0k}^{\cdot 1b} = C_{C\cdot 10}^{1} C_{AC\cdot ik}^{b0} + C_{C\cdot 10}^{0} C_{AC\cdot ik}^{b1}\\
&C_{AR\cdot\cdot 0i\cdot 0k}^{\cdot 1b} = C_{C\cdot 00}^{1} C_{AC\cdot ik}^{b0} + C_{C\cdot 00}^{0} C_{AC\cdot ik}^{b1}\\
&C_{AR\cdot\cdot 1i\cdot 1k}^{\cdot 1b} = C_{C\cdot 11}^{1} C_{AC\cdot ik}^{b0} + C_{C\cdot 11}^{0} C_{AC\cdot ik}^{b1}
\end{aligned}
\right. \tag{3.19}
$$

(3.16) следует из равенств (3.3), (3.19). □

Из равенства (3.16) следует, что

$$
\begin{cases}
C_{AR\cdot\cdot 0i\cdot 0k}^{\cdot 0b} = -C_{AR\cdot\cdot 1i\cdot 1k}^{\cdot 0b} = C_{AC\cdot ik}^{b0}\\
C_{AR\cdot\cdot 0i\cdot 1k}^{\cdot 0b} = C_{AR\cdot\cdot 1i\cdot 0k}^{\cdot 0b} = -C_{AC\cdot ik}^{b1}\\
C_{AR\cdot\cdot 0i\cdot 1k}^{\cdot 1b} = C_{AR\cdot\cdot 1i\cdot 0k}^{\cdot 1b} = C_{AC\cdot ik}^{b0}\\
C_{AR\cdot\cdot 0i\cdot 0k}^{\cdot 1b} = -C_{AR\cdot\cdot 1i\cdot 1k}^{\cdot 1b} = C_{AC\cdot ik}^{b1}
\end{cases}
\tag{3.20}
$$

## 4. Линейное отображение $C$-алгебры

**Определение 4.1.** Отображение $C$-векторных пространств

$$f : A_1 \to A_2$$

называется $C$-линейным, если

$$f(ca) = cf(a) \quad c \in C \quad a \in A \tag{4.1}$$

□

**Теорема 4.2** (Уравнения Коши-Римана)**.** *Матрица линейного отображения $C$-векторных пространств*

$$f : A_1 \to A_2$$
$$y^{ik} = f^{ik}_{jm} x^{jm}$$

*удовлетворяет соотношению*

$$
\begin{aligned}
f_{\cdot 0i}^{\cdot 0j} &= f_{\cdot 1i}^{\cdot 1j}\\
f_{\cdot 0i}^{\cdot 1j} &= -f_{\cdot 1i}^{\cdot 0j}
\end{aligned}
\tag{4.2}
$$

*Доказательство.* Из равенств (3.11), (3.14), (3.15) следует

$$
\begin{aligned}
(f(ca))^{0j} &= f^{0j}_{0i}(ca)^{0i} + f^{0j}_{1i}(ca)^{1i}\\
&= f^{0j}_{0i}(c^0 a^{0i} - c^1 a^{1i}) + f^{0j}_{1i}(c^0 a^{1i} + c^1 a^{0i})\\
&= f^{0j}_{0i} c^0 a^{0i} - f^{0j}_{0i} c^1 a^{1i} + f^{0j}_{1i} c^0 a^{1i} + f^{0j}_{1i} c^1 a^{0i}
\end{aligned}
\tag{4.3}
$$

$$
\begin{aligned}
(f(ca))^{1j} &= f^{1j}_{0i}(ca)^{0i} + f^{1j}_{1i}(ca)^{1i}\\
&= f^{1j}_{0i}(c^0 a^{0i} - c^1 a^{1i}) + f^{1j}_{1i}(c^0 a^{1i} + c^1 a^{0i})\\
&= f^{1j}_{0i} c^0 a^{0i} - f^{1j}_{0i} c^1 a^{1i} + f^{1j}_{1i} c^0 a^{1i} + f^{1j}_{1i} c^1 a^{0i}
\end{aligned}
\tag{4.4}
$$

Из равенств (4.1), (3.11), (3.12), (3.14), (3.15) следует

$$
\begin{aligned}
(cf(a))^{0j} &= c^0(f(a))^{0j} - c^1(f(a))^{1j}\\
&= c^0(f^{0j}_{0i} a^{0i} + f^{0j}_{1i} a^{1i}) - c^1(f^{1j}_{0i} a^{0i} + f^{1j}_{1i} a^{1i})\\
&= c^0 f^{0j}_{0i} a^{0i} + c^0 f^{0j}_{1i} a^{1i} - c^1 f^{1j}_{0i} a^{0i} - c^1 f^{1j}_{1i} a^{1i}
\end{aligned}
\tag{4.5}
$$

$$
\begin{aligned}
(cf(a))^{1j} &= c^0(f(a))^{1j} + c^1(f(a))^{0j}\\
&= c^0(f^{1j}_{0i} a^{0i} + f^{1j}_{1i} a^{1i}) + c^1(f^{0j}_{0i} a^{0i} + f^{0j}_{1i} a^{1i})\\
&= c^0 f^{1j}_{0i} a^{0i} + c^0 f^{1j}_{1i} a^{1i} + c^1 f^{0j}_{0i} a^{0i} + c^1 f^{0j}_{1i} a^{1i}
\end{aligned}
\tag{4.6}
$$

Из равенств (4.1), (4.3), (4.5) следует

$$
\begin{aligned}
& f^{0j}_{0i}c^0a^{0i} - f^{0j}_{0i}c^1a^{1i} + f^{0j}_{1i}c^0a^{1i} + f^{0j}_{1i}c^1a^{0i} \\
=& c^0f^{0j}_{0i}a^{0i} + c^0f^{0j}_{1i}a^{1i} - c^1f^{1j}_{0i}a^{0i} - c^1f^{1j}_{1i}a^{1i}
\end{aligned}
\tag{4.7}
$$

Из равенств (4.1), (4.4), (4.6) следует

$$
\begin{aligned}
& f^{1j}_{0i}c^0a^{0i} - f^{1j}_{0i}c^1a^{1i} + f^{1j}_{1i}c^0a^{1i} + f^{1j}_{1i}c^1a^{0i} \\
=& c^0f^{1j}_{0i}a^{0i} + c^0f^{1j}_{1i}a^{1i} + c^1f^{0j}_{0i}a^{0i} + c^1f^{0j}_{1i}a^{1i}
\end{aligned}
\tag{4.8}
$$

Так как $c^k$, $a^{ki}$ произвольны, то равенство (4.2) следует из равенств (4.7), (4.8). □

**Теорема 4.3.** *Пусть $\overline{\overline{e}}_{iC}$ - $C$-базис $C$-векторного пространства $A_i$, $i=1$, 2. Пусть $\overline{\overline{e}}_{iR}$ - соответствующий $R$-базис $C$-векторного пространства $A_i$. Матрица линейного отображения*

$$f:A_1\to A_2$$

*имеет вид*

$$
\begin{aligned}
x^{0j}_2 &= f^{0j}_{0i}x^{0i}_1 - f^{1j}_{0i}x^{1i}_1 \\
x^{1j}_2 &= f^{1j}_{0i}x^{0i}_1 + f^{0j}_{0i}x^{1i}_1
\end{aligned}
\tag{4.9}
$$

*относительно базисов $\overline{\overline{e}}_{1R}$, $\overline{\overline{e}}_{2R}$ и имеет вид*

$$x^j_2 = f^j_i x^i_1 \quad f^j_i = f^{0j}_{0i} + f^{1j}_{0i}i \tag{4.10}$$

*относительно базисов $\overline{\overline{e}}_{1C}$, $\overline{\overline{e}}_{2C}$*

*Доказательство.* Так как отображение $f$ является $R$-линейным отображением, то относительно базисов $\overline{\overline{e}}_{1R}$, $\overline{\overline{e}}_{2R}$ отображение $f$ имеет вид

$$
\begin{aligned}
x^{0j}_2 &= f^{0j}_{0i}x^{0i}_1 + f^{0j}_{1i}x^{1i}_1 \\
x^{1j}_2 &= f^{1j}_{0i}x^{0i}_1 + f^{1j}_{1i}x^{1i}_1
\end{aligned}
\tag{4.11}
$$

Равенство (4.9) следует из равенств (4.2), (4.11). Равенство (4.10) следует из равенства (4.9) и сравнения равенств

$$(f^{0j}_{0i}+f^{1j}_{0i}i)(x^{0i}_1+x^{1i}_1 i) = f^{0j}_{0i}x^{0i}_1 - f^{1j}_{0i}x^{1i}_1 + (f^{1j}_{0i}x^{0i}_1 + f^{0j}_{0i}x^{1i}_1)i$$

$$\begin{pmatrix} f^{0j}_{0i} & -f^{1j}_{0i} \\ f^{1j}_{0i} & f^{0j}_{0i} \end{pmatrix} \begin{pmatrix} x^{0i}_1 \\ x^{1i}_1 \end{pmatrix} = \begin{pmatrix} f^{0j}_{0i}x^{0i}_1 - f^{1j}_{0i}x^{1i}_1 \\ f^{1j}_{0i}x^{0i}_1 + f^{0j}_{0i}x^{1i}_1 \end{pmatrix}$$

□

Теорема 4.2 описывает структуру $C$-линейного отображения $C$-алгебры $A$. Однако существуют $R$-линейные отображения $C$-алгебры $A$, не являющиеся $C$-линейными.

## 5. Антилинейное отображение $C$-алгебры

**Определение 5.1.** Морфизм $C$-векторных пространств

(5.1) $$I : z \in C \to z^* \in C \quad f : A_1 \to A_2$$

называется антилинейным отображением.[1] □

**Теорема 5.2.** *Пусть отображение $C$-векторных пространств*

$$f : A_1 \to A_2$$

*является антилинейным отображением. Тогда*

(5.2) $$f(ca) = c^* f(a) \quad c \in C \quad a \in A$$

*Доказательство.* Согласно построениям, выполненным в разделе [4]-1.5, отображение $I$ является морфизмом $R$-алгебры $C$. Равенство (5.2) является следствием определения [5]-2.2.2, и равенства [5]-(2.2.4). □

Если $c \in R$, то

$$f(ca) = cf(a)$$

так как $c = c^*$. Следовательно, антилинейное отображение $C$-алгебры $A$ является $R$-линейным отображением.

**Теорема 5.3** (Уравнения Коши-Римана)**.** *Матрица антилинейного отображения $C$-векторных пространств*

$$f : A_1 \to A_2$$
$$y^{\boldsymbol{ik}} = f^{\boldsymbol{ik}}_{\boldsymbol{jm}} x^{\boldsymbol{jm}}$$

*удовлетворяет соотношению*

(5.3) $$\begin{aligned} f^{\boldsymbol{0j}}_{\boldsymbol{0i}} &= -f^{\boldsymbol{1j}}_{\boldsymbol{1i}} \\ f^{\boldsymbol{0j}}_{\boldsymbol{1i}} &= \phantom{-}f^{\boldsymbol{1j}}_{\boldsymbol{0i}} \end{aligned}$$

*Доказательство.* Из равенств (3.11), (3.14), (3.15) следует

(5.4) $$\begin{aligned} (f(ca))^{\boldsymbol{0j}} &= f^{\boldsymbol{0j}}_{\boldsymbol{0i}}(ca)^{\boldsymbol{0i}} + f^{\boldsymbol{0j}}_{\boldsymbol{1i}}(ca)^{\boldsymbol{1i}} \\ &= f^{\boldsymbol{0j}}_{\boldsymbol{0i}}(c^{\boldsymbol{0}}a^{\boldsymbol{0i}} - c^{\boldsymbol{1}}a^{\boldsymbol{1i}}) + f^{\boldsymbol{0j}}_{\boldsymbol{1i}}(c^{\boldsymbol{0}}a^{\boldsymbol{1i}} + c^{\boldsymbol{1}}a^{\boldsymbol{0i}}) \\ &= f^{\boldsymbol{0j}}_{\boldsymbol{0i}}c^{\boldsymbol{0}}a^{\boldsymbol{0i}} - f^{\boldsymbol{0j}}_{\boldsymbol{0i}}c^{\boldsymbol{1}}a^{\boldsymbol{1i}} + f^{\boldsymbol{0j}}_{\boldsymbol{1i}}c^{\boldsymbol{0}}a^{\boldsymbol{1i}} + f^{\boldsymbol{0j}}_{\boldsymbol{1i}}c^{\boldsymbol{1}}a^{\boldsymbol{0i}} \end{aligned}$$

(5.5) $$\begin{aligned} (f(ca))^{\boldsymbol{1j}} &= f^{\boldsymbol{1j}}_{\boldsymbol{0i}}(ca)^{\boldsymbol{0i}} + f^{\boldsymbol{1j}}_{\boldsymbol{1i}}(ca)^{\boldsymbol{1i}} \\ &= f^{\boldsymbol{1j}}_{\boldsymbol{0i}}(c^{\boldsymbol{0}}a^{\boldsymbol{0i}} - c^{\boldsymbol{1}}a^{\boldsymbol{1i}}) + f^{\boldsymbol{1j}}_{\boldsymbol{1i}}(c^{\boldsymbol{0}}a^{\boldsymbol{1i}} + c^{\boldsymbol{1}}a^{\boldsymbol{0i}}) \\ &= f^{\boldsymbol{1j}}_{\boldsymbol{0i}}c^{\boldsymbol{0}}a^{\boldsymbol{0i}} - f^{\boldsymbol{1j}}_{\boldsymbol{0i}}c^{\boldsymbol{1}}a^{\boldsymbol{1i}} + f^{\boldsymbol{1j}}_{\boldsymbol{1i}}c^{\boldsymbol{0}}a^{\boldsymbol{1i}} + f^{\boldsymbol{1j}}_{\boldsymbol{1i}}c^{\boldsymbol{1}}a^{\boldsymbol{0i}} \end{aligned}$$

[1] Аналогичное определение дано в [7], с. 234.

Из равенств (5.2), (3.11), (3.13), (3.14), (3.15) следует

$$\begin{aligned}(c^*f(a))^{0j} &= c^0(f(a))^{0j}+c^1(f(a))^{1j}\\ &= c^0(f^{0j}_{0i}a^{0i}+f^{0j}_{1i}a^{1i})+c^1(f^{1j}_{0i}a^{0i}+f^{1j}_{1i}a^{1i})\\ &= c^0f^{0j}_{0i}a^{0i}+c^0f^{0j}_{1i}a^{1i}+c^1f^{1j}_{0i}a^{0i}+c^1f^{1j}_{1i}a^{1i}\end{aligned} \tag{5.6}$$

$$\begin{aligned}(c^*f(a))^{1j} &= c^0(f(a))^{1j}-c^1(f(a))^{0j}\\ &= c^0(f^{1j}_{0i}a^{0i}+f^{1j}_{1i}a^{1i})-c^1(f^{0j}_{0i}a^{0i}+f^{0j}_{1i}a^{1i})\\ &= c^0f^{1j}_{0i}a^{0i}+c^0f^{1j}_{1i}a^{1i}-c^1f^{0j}_{0i}a^{0i}-c^1f^{0j}_{1i}a^{1i}\end{aligned} \tag{5.7}$$

Из равенств (5.2), (5.4), (5.6) следует

$$\begin{aligned}&f^{0j}_{0i}c^0a^{0i}-f^{0j}_{0i}c^1a^{1i}+f^{0j}_{1i}c^0a^{1i}+f^{0j}_{1i}c^1a^{0i}\\ =&c^0f^{0j}_{0i}a^{0i}+c^0f^{0j}_{1i}a^{1i}+c^1f^{1j}_{0i}a^{0i}+c^1f^{1j}_{1i}a^{1i}\end{aligned} \tag{5.8}$$

Из равенств (5.2), (5.5), (5.7) следует

$$\begin{aligned}&f^{1j}_{0i}c^0a^{0i}-f^{1j}_{0i}c^1a^{1i}+f^{1j}_{1i}c^0a^{1i}+f^{1j}_{1i}c^1a^{0i}\\ =&c^0f^{1j}_{0i}a^{0i}+c^0f^{1j}_{1i}a^{1i}-c^1f^{0j}_{0i}a^{0i}-c^1f^{0j}_{1i}a^{1i}\end{aligned} \tag{5.9}$$

Так как $c^k$, $a^{ki}$ произвольны, то равенство (5.3) следует из равенств (5.8), (5.9). □

**Теорема 5.4.** *Пусть $\overline{\overline{e}}_{iC}$ - $C$-базис $C$-векторного пространства $A_i$, $i=1$, 2. Пусть $\overline{\overline{e}}_{iR}$ - соответствующий $R$-базис $C$-векторного пространства $A_i$. Матрица антилинейного отображения*

$$f:A_1\to A_2$$

*имеет вид*

$$\begin{aligned}x_2^{0j}&=f^{0j}_{0i}x_1^{0i}+f^{1j}_{0i}x_1^{1i}\\ x_2^{1j}&=f^{1j}_{0i}x_1^{0i}-f^{0j}_{0i}x_1^{1i}\end{aligned} \tag{5.10}$$

*относительно базисов $\overline{\overline{e}}_{1R}$, $\overline{\overline{e}}_{2R}$ и имеет вид*

$$x_2^j=f^j_iI\circ x_1^i\quad I\circ z=z^*\quad f^j_i=f^{0j}_{0i}+f^{1j}_{0i}i \tag{5.11}$$

*относительно базисов $\overline{\overline{e}}_{1C}$, $\overline{\overline{e}}_{2C}$*

*Доказательство.* Так как отображение $f$ является $R$-линейным отображением, то относительно базисов $\overline{\overline{e}}_{1R}$, $\overline{\overline{e}}_{2R}$ отображение $f$ имеет вид

$$\begin{aligned}x_2^{0j}&=f^{0j}_{0i}x_1^{0i}+f^{0j}_{1i}x_1^{1i}\\ x_2^{1j}&=f^{1j}_{0i}x_1^{0i}+f^{1j}_{1i}x_1^{1i}\end{aligned} \tag{5.12}$$

Равенство (5.10) следует из равенств (5.3), (5.12). Равенство (5.11) следует из равенства (5.10) и сравнения равенств

$$\begin{aligned}(f^{0j}_{0i}+f^{1j}_{0i}i)I\circ(x_1^{0i}+x_1^{1i}i)&=(f^{0j}_{0i}+f^{1j}_{0i}i)(x_1^{0i}-x_1^{1i}i)\\ &=f^{0j}_{0i}x_1^{0i}+f^{1j}_{0i}x_1^{1i}+(f^{1j}_{0i}x_1^{0i}-f^{0j}_{0i}x_1^{1i})i\end{aligned}$$

$$\begin{pmatrix} f^{0j}_{0i} & f^{1j}_{0i} \\ f^{1j}_{0i} & -f^{0j}_{0i} \end{pmatrix} \begin{pmatrix} x^{0i}_1 \\ x^{1i}_1 \end{pmatrix} = \begin{pmatrix} f^{0j}_{0i} x^{0i}_1 + f^{1j}_{0i} x^{1i}_1 \\ f^{1j}_{0i} x^{0i}_1 - f^{0j}_{0i} x^{1i}_1 \end{pmatrix}$$

□

## 6. Морфизм алгебры

В этом разделе мы рассмотрим алгебру $A$ над кольцом $D$. Согласно построениям, выполненным в разделах [5]-4.4.2, [5]-4.4.3, диаграмма представлений $D$-алгебры имеет вид

(6.1)
$$\begin{array}{ccccc} D & \xrightarrow[*]{f_{1,2}} & A & \xrightarrow[*]{f_{2,3}} & A \\ & & & & \uparrow{\scriptstyle *\, f_{1,2}} \\ & & & & D \end{array} \qquad \begin{array}{r} f_{1,2}(d) : \overline{v} \to d\,\overline{v} \\ f_{2,3}(\overline{v}) : \overline{w} \to \overline{C}(\overline{v}, \overline{w}) \\ \overline{C} \in \mathcal{L}(A^2; A) \end{array}$$

На диаграмме представлений (6.1), $D$ - кольцо, $A$ - абелевая группа. Мы сперва рассматриваем вертикальное представление, а затем горизонтальные. Мы полагаем, что $A$ - свободный $D\star$-модуль ([1], с. 103). Следовательно, векторы базиса $\overline{\overline{e}}$ $D\star$-модуля $A$ линейно независимы. Билинейное отображение

$$\overline{C} : A \times A \to A$$

имеет следующее разложение относительно базиса $\overline{\overline{e}}$

(6.2)
$$\overline{C}(\overline{v}, \overline{w}) = C^{k}_{ij} v^{i} w^{j} \overline{e}_{k}$$

Согласно теореме [5]-4.4.1, билинейное отображение $\overline{C}$ порождает произведение в $D\star$-алгебре $A$

(6.3)
$$\overline{C}(\overline{v}, \overline{w}) = \overline{v}\,\overline{w}$$

и координаты билинейного отображения $\overline{C}$ относительно базиса $\overline{\overline{e}}$ являются структурными константами этого произведения.

**Теорема 6.1.** *Пусть*

(6.4)
$$\begin{array}{ccccc} D_1 & \xrightarrow[*]{f_{1\cdot1,2}} & A_1 & \xrightarrow[*]{f_{1\cdot2,3}} & A_1 \\ & & & & \uparrow{\scriptstyle *\, f_{1\cdot1,2}} \\ & & & & D_1 \end{array} \qquad \begin{array}{r} f_{1\cdot1,2}(d) : \overline{v} \to d\,\overline{v} \\ f_{1\cdot2,3}(\overline{v}) : \overline{w} \to \overline{C}_1(\overline{v}, \overline{w}) \\ \overline{C}_1 \in \mathcal{L}(A_1^2; A_1) \end{array}$$

*диаграмма представлений, описывающая $D_1\star$-алгебру $A_1$. Пусть*

(6.5)
$$\begin{array}{ccccc} D_2 & \xrightarrow[*]{f_{2\cdot1,2}} & A_2 & \xrightarrow[*]{f_{2\cdot2,3}} & A_2 \\ & & & & \uparrow{\scriptstyle *\, f_{2\cdot1,2}} \\ & & & & D_2 \end{array} \qquad \begin{array}{r} f_{2\cdot1,2}(d) : \overline{v} \to d\,\overline{v} \\ f_{2\cdot2,3}(\overline{v}) : \overline{w} \to \overline{C}_2(\overline{v}, \overline{w}) \\ \overline{C}_2 \in \mathcal{L}(A_2^2; A_2) \end{array}$$

*диаграмма представлений, описывающая $D_2\star$-алгебру $A_2$. Морфизм*[2] *$D_1\star$-алгебры $A_1$ в $D_2\star$-алгебру $A_2$*

$$r_1 : D_1 \to D_2 \quad \overline{r}_2 : A_1 \to A_2$$

*является $D\star$-линейным отображением $D_1\star$-алгебры $A_1$ в $D_2\star$-алгебру $A_2$ таким, что*

$$\overline{r}_2(\overline{a}\overline{b}) = \overline{r}_2(\overline{a})\overline{r}_2(\overline{b}) \tag{6.6}$$

*Доказательство.* Согласно равенству [5]-(4.2.3), морфизм $(r_1, \overline{r}_2)$ представления $f_{1,2}$ удовлетворяет равенству

$$\begin{aligned} \overline{r}_2(f_{1\cdot 1,2}(d)(\overline{a})) &= f_{2\cdot 1,2}(r_1(d))(\overline{r}_2(\overline{a})) \\ \overline{r}_2(d\,\overline{a}) &= r_1(d)\overline{r}_2(\overline{a}) \end{aligned} \tag{6.7}$$

Следовательно, отображение $(r_1, \overline{r}_2)$ является $D\star$-линейным отображением.

Согласно равенству [5]-(4.2.3), морфизм $(\overline{r}_2, \overline{r}_2)$ представления $f_{2,3}$ удовлетворяет равенству[3]

$$\overline{r}_2(f_{1\cdot 2,3}(\overline{a}_2)(\overline{a}_3)) = f_{2\cdot 2,3}(\overline{r}_2(\overline{a}_2))(\overline{r}_2(\overline{a}_3)) \tag{6.8}$$

Из равенств (6.8), (6.4), (6.5), следует

$$\overline{r}_2(\overline{C}_1(\overline{v}, \overline{w})) = \overline{C}_2(\overline{r}_2(\overline{v}), \overline{r}_2(\overline{w})) \tag{6.9}$$

Равенство (6.6) следует из равенств (6.9), (6.3). □

**Определение 6.2.** Морфизм представлений $D_1\star$-алгебры $A_1$ в $D_2\star$-алгебру $A_2$ называется **$D\star$-линейным гомоморфизмом** $D_1\star$-алгебры $A_1$ в $D_2\star$-алгебру $A_2$. □

**Теорема 6.3.** *Пусть $\overline{\overline{e}}_1$ - базис $D_1\star$-алгебры $A_1$. Пусть $\overline{\overline{e}}_2$ - базис $D_2\star$-алгебры $A_2$. Тогда линейный гомоморфизм*[4] *$(r_1, \overline{r}_2)$ $D_1\star$-алгебры $A_1$ в $D_2\star$-алгебру $A_2$ имеет представление*

$$\begin{aligned} \overline{b} &= r_1(a)^*{}_*r_2{}^*{}_*\overline{e}_2 = r_1(a^j) r_{2\cdot}{}^i_j\, \overline{e}_{2\cdot i} \\ b &= r_1(a)^*{}_*r_2 \end{aligned} \tag{6.10}$$

*относительно заданных базисов. Здесь*

- *$a$ - координатная матрица вектора $\overline{a}$ относительно базиса $\overline{\overline{e}}_1$.*
- *$b$ - координатная матрица вектора*

$$\overline{b} = \overline{r}_2(\overline{a})$$

*относительно базиса $\overline{\overline{e}}_2$.*

---

[2]В замечании [2]-4.4.5, я показал, что выбор отображения $r_1$ при изучении линейных отображений может быть ограничен отображением

$$r_1 : D \to D \quad r_1(d) = d$$

Однако в разделе 10 я рассматриваю антилинейные отображения алгебры. Поэтому я записал теоремы в этом разделе таким образом, чтобы их можно было сопоставить с теоремами в разделе 10.

[3]Так как в диаграммах представлений (6.4), (6.5), носители $\Omega_2$-алгебры и $\Omega_3$-алгебры совпадают, то также совпадают морфизмы представлений на уровнях 2 и 3.

[4]Эта теорема аналогична теореме [2]-4.4.3.

- $r_2$ - *координатная матрица множества векторов* $(\overline{r}_2(\overline{e}_{1\cdot i}))$ *относительно базиса* $\overline{\overline{e}}_2$. *Мы будем называть матрицу* $r_2$ **матрицей линейного гомоморфизма** *относительно базисов* $\overline{\overline{e}}_1$ *и* $\overline{\overline{e}}_2$.

*Доказательство.* Вектор $\overline{a}\in A_1$ имеет разложение

$$\overline{a}=a^*{}_*\overline{e}_1$$

относительно базиса $\overline{\overline{e}}_1$. Вектор $\overline{b}\in A_2$ имеет разложение

$$\overline{b}=b^*{}_*\overline{e}_2 \tag{6.11}$$

относительно базиса $\overline{\overline{e}}_2$.

Так как $(r_1,\overline{r}_2)$ - линейный гомоморфизм, то на основании (6.7) следует, что

$$\overline{b}=\overline{r}_2(\overline{a})=\overline{r}_2(a^*{}_*\overline{e}_1)=r_1(a)^*{}_*\overline{r}_2(\overline{e}_1) \tag{6.12}$$

где

$$r_1(a)=\begin{pmatrix}r_1(a^1)\\ ...\\ r_1(a^n)\end{pmatrix}$$

$\overline{r}_2(\overline{e}_{1\cdot i})$ также вектор $D\star$-модуля $A_2$ и имеет разложение

$$\overline{r}_2(\overline{e}_{1\cdot i})=r_{2\cdot i}{}^*{}_*\overline{e}_2=r_{2\cdot}{}^j_i\,\overline{e}_j \tag{6.13}$$

относительно базиса $\overline{\overline{e}}_2$. Комбинируя (6.12) и (6.13), мы получаем

$$\overline{b}=r_1(a)^*{}_*r_2{}^*{}_*\overline{e}_2 \tag{6.14}$$

(6.10) следует из сравнения (6.11) и (6.14) и теоремы [2]-4.3.3. □

**Теорема 6.4.** *Пусть* $\overline{\overline{e}}_1$ *- базис* $D_1\star$*-алгебры* $A_1$. *Пусть* $\overline{\overline{e}}_2$ *- базис* $D_2\star$*-алгебры* $A_2$. *Если отображение* $r_1$ *является инъекцией, то матрица линейного гомоморфизма и структурные константы связаны соотношением*

$$r_1(C_{1\cdot}{}^k_{ij})r_{2\cdot}{}^l_k=r_{2\cdot}{}^p_i r_{2\cdot}{}^q_j C_{2\cdot}{}^l_{pq} \tag{6.15}$$

*Доказательство.* Пусть

$$\overline{a},\overline{b}\in A_1\quad \overline{a}=a^*{}_*\overline{e}_1\quad \overline{b}=b^*{}_*\overline{e}_1$$

Из равенств (6.2), (6.3), (6.4), следует

$$\overline{a}\overline{b}=a^ib^jC_{1\cdot}{}^k_{ij}\overline{e}_{1\cdot k} \tag{6.16}$$

Из равенств (6.7), (6.16), следует

$$\overline{r}_2(\overline{a}\overline{b})=r_1(C_{1\cdot}{}^k_{ij}a^ib^j)\overline{r}_2(\overline{e}_{1\cdot k}) \tag{6.17}$$

Поскольку отображение $r_1$ является гомоморфизмом колец, то из равенства (6.17), следует

$$\overline{r}_2(\overline{a}\overline{b})=r_1(C_{1\cdot}{}^k_{ij})r_1(a^i)r_1(b^j)\overline{r}_2(\overline{e}_{1\cdot k}) \tag{6.18}$$

Из теоремы 6.3 и равенства (6.18), следует

$$\overline{r}_2(\overline{a}\overline{b})=r_1(C_{1\cdot}{}^k_{ij})r_1(a^i)r_1(b^j)r_{2\cdot}{}^l_k\,\overline{e}_{2\cdot l} \tag{6.19}$$

Из равенства (6.6) и теоремы 6.3, следует

$$\overline{r}_2(\overline{a}\overline{b})=\overline{r}_2(\overline{a})\overline{r}_2(\overline{b})=r_1(a^i)r_{2\cdot}{}^p_i\,\overline{e}_{2\cdot p}r_1(b^j)r_{2\cdot}{}^q_j\,\overline{e}_{2\cdot q} \tag{6.20}$$

Из равенств (6.2), (6.3), (6.5), (6.20), следует

(6.21) $$\overline{r}_2(\overline{a}\overline{b}) = r_1(a^i)r_{2\cdot}{}^p_i r_1(b^j)r_{2\cdot}{}^q_j C_{2\cdot}{}^l_{pq}\,\overline{e}_{2\cdot l}$$

Из равенств (6.19), (6.21), следует

(6.22) $$r_1(C_{1\cdot}{}^k_{ij})r_1(a^i)r_1(b^j)r_{2\cdot}{}^l_k\,\overline{e}_{2\cdot l} = r_1(a^i)r_{2\cdot}{}^p_i r_1(b^j)r_{2\cdot}{}^q_j C_{2\cdot}{}^l_{pq}\,\overline{e}_{2\cdot l}$$

Равенство (6.15) следует из равенства (6.22), так как векторы базиса $\overline{\overline{e}}_2$ линейно независимы, и $a^i$, $b^i$ (а следовательно, $r_1(a^i)$, $r_1(b^i)$) - произвольные величины. □

## 7. Линейный автоморфизм алгебры кватернионов

Определение координат линейного автоморфизма - задача непростая. В этом разделе мы рассмотрим пример нетривиального линейного автоморфизма алгебры кватернионов.

**Теорема 7.1.** *Координаты линейного автоморфизма алгебры кватернионов удовлетворяют системе уравнений*

(7.1) $$\begin{cases} r^1_1 = r^2_2r^3_3 - r^3_2r^2_3 & r^1_2 = r^2_3r^3_1 - r^3_3r^2_1 & r^1_3 = r^2_1r^3_2 - r^3_1r^2_2 \\ r^2_1 = r^3_2r^1_3 - r^1_2r^3_3 & r^2_2 = r^3_3r^1_1 - r^1_3r^3_1 & r^2_3 = r^3_1r^1_2 - r^1_1r^3_2 \\ r^3_1 = r^1_2r^2_3 - r^2_2r^1_3 & r^3_2 = r^1_3r^2_1 - r^2_3r^1_1 & r^3_3 = r^1_1r^2_2 - r^2_1r^1_2 \end{cases}$$

*Доказательство.* Согласно теоремам [4]-3.3.1, 6.4, линейный автоморфизм алгебры кватернионов удовлетворяет уравнениям

(7.2) $$\begin{array}{llll} r^l_0 = r^p_0r^q_0C^l_{pq} & r^l_1 = r^p_0r^q_1C^l_{pq} & r^l_2 = r^p_0r^q_2C^l_{pq} & r^l_3 = r^p_0r^q_3C^l_{pq} \\ r^l_1 = r^p_1r^q_0C^l_{pq} & -r^l_0 = r^p_1r^q_1C^l_{pq} & r^l_3 = r^p_1r^q_2C^l_{pq} & -r^l_2 = r^p_1r^q_3C^l_{pq} \\ r^l_2 = r^p_2r^q_0C^l_{pq} & -r^l_3 = r^p_2r^q_1C^l_{pq} & -r^l_0 = r^p_2r^q_2C^l_{pq} & r^l_0 = r^p_2r^q_3C^l_{pq} \\ r^l_3 = r^p_3r^q_0C^l_{pq} & r^l_2 = r^p_3r^q_1C^l_{pq} & -r^l_1 = r^p_3r^q_2C^l_{pq} & -r^l_0 = r^p_3r^q_3C^l_{pq} \end{array}$$

Из равенства (7.2) следует

(7.3) $$\begin{array}{lll} r^l_1 = r^p_0r^q_1C^l_{pq} = r^p_2r^q_3C^l_{pq} & r^p_0r^q_1C^l_{pq} = r^p_0r^q_1C^l_{qp} & r^p_2r^q_3C^l_{pq} = -r^p_2r^q_3C^l_{qp} \\ r^l_2 = r^p_0r^q_2C^l_{pq} = r^p_3r^q_1C^l_{pq} & r^p_0r^q_2C^l_{pq} = r^p_0r^q_2C^l_{qp} & r^p_3r^q_1C^l_{pq} = -r^p_1r^q_3C^l_{pq} \\ r^l_3 = r^p_0r^q_3C^l_{pq} = r^p_1r^q_2C^l_{pq} & r^p_0r^q_3C^l_{pq} = r^p_0r^q_3C^l_{qp} & r^p_1r^q_2C^l_{pq} = -r^p_1r^q_2C^l_{qp} \end{array}$$

(7.4) $$r^l_0 = r^p_0r^q_0C^l_{pq} = -r^p_1r^q_1C^l_{pq} = -r^p_2r^q_2C^l_{pq} = -r^p_3r^q_3C^l_{pq}$$

Если $l = 0$, то из равенства

$$C^0_{pq} = C^0_{qp}$$

следует

(7.5) $$r^p_ir^q_jC^0_{pq} = r^p_ir^q_jC^0_{qp}$$

Из равенства (7.3) для $l = 0$ и равенства (7.5), следует

(7.6) $$r^0_1 = r^0_2 = r^0_3 = 0$$

Если $l=1, 2, 3,$ то равенство (7.3) можно записать в виде

(7.7)
$$
\left\{
\begin{aligned}
r_i^l &= r_0^l r_i^0 C_{l0}^l + r_0^0 r_i^l C_{0l}^l + r_0^a r_i^b C_{ab}^l + r_0^b r_i^a C_{ba}^l \\
&\quad r_0^l r_i^0 C_{l0}^l + r_0^0 r_i^l C_{0l}^l + r_0^a r_i^b C_{ab}^l + r_0^b r_i^a C_{ba}^l \\
&= r_0^l r_i^0 C_{0l}^l + r_0^0 r_i^l C_{l0}^l + r_0^a r_i^b C_{ba}^l + r_0^b r_i^a C_{ab}^l \\
&\quad i = 1, 2, 3 \\
r_i^l &= r_k^0 r_j^l C_{0l}^l + r_k^l r_j^0 C_{l0}^l + r_k^a r_j^b C_{ab}^l + r_k^b r_j^a C_{ba}^l \\
&\quad r_k^0 r_j^l C_{0l}^l + r_k^l r_j^0 C_{l0}^l + r_k^a r_j^b C_{ab}^l + r_k^b r_j^a C_{ba}^l \\
&= -r_k^0 r_j^l C_{l0}^l - r_k^l r_j^0 C_{0l}^l - r_k^a r_j^b C_{ba}^l - r_k^b r_j^a C_{ab}^l \\
&\qquad i = 1 \quad k = 2 \quad j = 3 \\
&\qquad i = 2 \quad k = 3 \quad j = 1 \\
&\qquad i = 3 \quad k = 1 \quad j = 2 \\
&\quad 0 < a < b \quad a \ne l \quad b \ne l
\end{aligned}
\right.
$$

Из равенств (7.7), (7.6) и равенств

(7.8)
$$
\begin{aligned}
C_{0l}^l &= C_{l0}^l = 1 \\
C_{ab}^l &= -C_{ba}^l
\end{aligned}
$$

следует

(7.9)
$$
\left\{
\begin{aligned}
r_i^l &= r_0^0 r_i^l + r_0^a r_i^b C_{ab}^l - r_0^b r_i^a C_{ab}^l \\
&\quad r_0^0 r_i^l + r_0^a r_i^b C_{ab}^l - r_0^b r_i^a C_{ab}^l \\
&= r_0^0 r_i^l - r_0^a r_i^b C_{ab}^l + r_0^b r_i^a C_{ab}^l \\
&\quad i = 1, 2, 3 \\
r_i^l &= r_k^a r_j^b C_{ab}^l - r_k^b r_j^a C_{ab}^l \\
&\quad r_k^a r_j^b C_{ab}^l - r_k^b r_j^a C_{ab}^l \\
&= r_k^a r_j^b C_{ab}^l - r_k^b r_j^a C_{ab}^l \\
&\qquad i = 1 \quad k = 2 \quad j = 3 \\
&\qquad i = 2 \quad k = 3 \quad j = 1 \\
&\qquad i = 3 \quad k = 1 \quad j = 2 \\
&\quad 0 < a < b \quad a \ne l \quad b \ne l
\end{aligned}
\right.
$$

Из равенств (7.9) следует

$$(7.10)\quad \begin{cases} r^{l}_{i} = r^{0}_{0} r^{l}_{i} \\ \quad r^{a}_{0} r^{b}_{i} - r^{b}_{0} r^{a}_{i} = 0 \\ \quad \boldsymbol{i = 1, 2, 3} \\ r^{l}_{i} = r^{a}_{k} r^{b}_{j} C^{l}_{ab} - r^{b}_{k} r^{a}_{j} C^{l}_{ab} \\ \qquad \boldsymbol{i = 1} \quad \boldsymbol{k = 2} \quad \boldsymbol{j = 3} \\ \qquad \boldsymbol{i = 2} \quad \boldsymbol{k = 3} \quad \boldsymbol{j = 1} \\ \qquad \boldsymbol{i = 3} \quad \boldsymbol{k = 1} \quad \boldsymbol{j = 2} \\ \quad \boldsymbol{0 < a < b} \quad \boldsymbol{a \neq l} \quad \boldsymbol{b \neq l} \end{cases}$$

Из равенства (7.10) следует

$$(7.11)\qquad r^{0}_{0} = 1$$

Из равенства (7.4) для $\boldsymbol{l = 0}$ следует

$$(7.12)\qquad \begin{aligned} r^{0}_{0} &= r^{0}_{0} r^{0}_{0} - r^{1}_{0} r^{1}_{0} - r^{2}_{0} r^{2}_{0} - r^{3}_{0} r^{3}_{0} \\ &= -r^{0}_{i} r^{0}_{i} + r^{1}_{i} r^{1}_{i} + r^{2}_{i} r^{2}_{i} + r^{3}_{i} r^{3}_{i} \\ \boldsymbol{i} &= \boldsymbol{1, 2, 3} \end{aligned}$$

Из равенств (7.6), (7.10), (7.12), следует

$$(7.13)\qquad \begin{aligned} 0 &= r^{1}_{0} r^{1}_{0} + r^{2}_{0} r^{2}_{0} + r^{3}_{0} r^{3}_{0} \\ 1 &= r^{1}_{i} r^{1}_{i} + r^{2}_{i} r^{2}_{i} + r^{3}_{i} r^{3}_{i} \\ \boldsymbol{i} &= \boldsymbol{1, 2, 3} \end{aligned}$$

Из равенств (7.13) следует[5]

$$(7.14)\qquad r^{1}_{0} = r^{2}_{0} = r^{3}_{0} = 0$$

Из равенства (7.4) для $\boldsymbol{l > 0}$ следует

$$(7.15)\qquad \begin{aligned} r^{l}_{0} &= r^{l}_{0} r^{0}_{0} C^{l}_{l0} + r^{0}_{0} r^{l}_{0} C^{l}_{0l} + r^{a}_{0} r^{b}_{0} C^{l}_{ab} + r^{b}_{0} r^{a}_{0} C^{l}_{ba} \\ &= -r^{l}_{i} r^{0}_{i} C^{l}_{l0} - r^{0}_{i} r^{l}_{i} C^{l}_{0l} - r^{a}_{i} r^{b}_{i} C^{l}_{ab} - r^{b}_{i} r^{a}_{i} C^{l}_{ba} \\ & \boldsymbol{i > 0} \\ & \boldsymbol{l > 0} \quad \boldsymbol{0 < a < b} \quad \boldsymbol{a \neq l} \quad \boldsymbol{b \neq l} \end{aligned}$$

---

[5]Мы здесь опираемся на то, что алгебра кватернионов определена над полем действительных чисел. Если рассматривать алгебру кватернионов над полем комплексных чисел, то уравнение (7.13) определяет конус в комплексном пространстве. Соответственно, у нас шире выбор координат линейного автоморфизма.

Равенства (7.15) тождественно верны в силу равенств (7.6), (7.14), (7.8). Из равенств (7.14), (7.10), следует

$$\begin{cases} r_i^l = r_k^a r_j^b C_{ab}^l - r_k^b r_j^a C_{ab}^l \\ i=1 \quad k=2 \quad j=3 \\ i=2 \quad k=3 \quad j=1 \\ i=3 \quad k=1 \quad j=2 \\ l>0 \quad 0<a<b \quad a\neq l \quad b\neq l \end{cases} \tag{7.16}$$

Равенства (7.1) следуют из равенств (7.16). □

**Пример 7.2.** Очевидно, координаты

$$r_j^i = \delta_j^i$$

удовлетворяют уравнению (7.1). Мы можем убедиться непосредственной проверкой, что координаты отображения

$$r_0^0 = 1 \quad r_2^1 = 1 \quad r_3^2 = 1 \quad r_1^3 = 1$$

также удовлетворяют уравнению (7.1). Матрица координат этого отображения имеет вид

$$r = \begin{pmatrix} 1 & 0 & 0 & \\ 0 & 0 & 1 & 0 \\ 0 & 0 & 0 & 1 \\ 0 & 1 & 0 & 0 \end{pmatrix}$$

Согласно теореме [4]-3.3.4, стандартные компоненты отображения $r$ имеют вид

$$\begin{array}{cccc} r^{00} = \frac{1}{4} & r^{11} = -\frac{1}{4} & r^{22} = -\frac{1}{4} & r^{33} = -\frac{1}{4} \\ r^{10} = -\frac{1}{4} & r^{01} = \frac{1}{4} & r^{32} = -\frac{1}{4} & r^{23} = -\frac{1}{4} \\ r^{20} = -\frac{1}{4} & r^{31} = -\frac{1}{4} & r^{02} = \frac{1}{4} & r^{13} = -\frac{1}{4} \\ r^{30} = -\frac{1}{4} & r^{21} = -\frac{1}{4} & r^{12} = -\frac{1}{4} & r^{03} = \frac{1}{4} \end{array}$$

Следовательно, отображение $r$ имеет вид

$$r(a) = a^0 + a^2 i + a^3 j + a^1 k$$

$$\begin{aligned} r(a) = \frac{1}{4}(&a - iai - jaj - kak - ia + ai - kaj - jak \\ &-ja - kai + aj - iak - ka - jai - iaj + ak) \end{aligned}$$

□

## 8. Кольцо с сопряжением

Пусть существует коммутативное подкольцо $F$ кольца $D$ такое, что $F \neq D$ и кольцо $D$ является свободным конечно мерным модулем над кольцом $F$. Пусть $\overline{\overline{e}}$ - базис свободного модуля $D$ над кольцом $F$. Мы будем полагать $\overline{e}_0 = 1$.

**Теорема 8.1.** *Кольцо $D$ является $F$-алгеброй.*

*Доказательство.* Рассмотрим отображение

$$f : D \times D \to D \quad a, b \in D \to ab \in D \tag{8.1}$$

Поскольку произведение дистрибутивно по отношению к сложению, то

$$\begin{aligned} f(a_1 + a_2, b) &= (a_1 + a_2)b = a_1 b + a_2 b = f(a_1, b) + f(a_2, b) \\ f(a, b_1 + b_2) &= a(b_1 + b_2) = ab_1 + ab_2 = f(a, b_1) + f(a, b_2) \end{aligned} \tag{8.2}$$

Если $a \in F$, то

$$\begin{aligned} f(ab, c) &= (ab)c = a(bc) = af(b, c) \\ f(b, ac) &= b(ac) = (ba)c = (ab)c = a(bc) = af(b, c) \end{aligned} \tag{8.3}$$

Из равенств (8.2), (8.3), следует, что отображение $f$ является билинейным отображением над кольцом $F$. Согласно определению [4]-2.2.1, кольцо $D$ является $F$-алгеброй. □

**Следствие 8.2.** *Кольцо $F$ является подалгеброй $F$-алгебры $D$.* □

Рассмотрим отображения

$$\begin{aligned} \mathrm{Re} &: D \to D \\ \mathrm{Im} &: D \to D \end{aligned}$$

определенные равенством

$$\mathrm{Re}\, d = d^0 \quad \mathrm{Im}\, d = d - d^0 \quad d \in D \quad d = d^i \overline{e}_i \tag{8.4}$$

Выражение $\mathrm{Re}\, d$ называется **скаляром элемента** $d$. Выражение $\mathrm{Im}\, d$ называется **вектором элемента** $d$.

Согласно (8.4)

$$F = \{d \in D : \mathrm{Re}\, d = d\}$$

Мы будем пользоваться записью $\mathrm{Re}\, D$ для обозначения **алгебры скаляров кольца** $D$.

**Теорема 8.3.** *Множество*

$$\mathrm{Im}\, D = \{d \in D : \mathrm{Re}\, d = 0\} \tag{8.5}$$

*является* $(\mathrm{Re}\, D)$*-модулем, который мы называем* **модуль векторов кольца** $D$*.*

$$D = \mathrm{Re}\, D \oplus \mathrm{Im}\, D \tag{8.6}$$

*Доказательство.* Пусть $c$, $d \in \mathrm{Im}\, D$. Тогда $c_0 = d_0 = 0$. Следовательно,

$$(c+d)_0 = c_0 + d_0 = 0$$

Если $a \in \mathrm{Re}\, D$, то

$$(ad)_0 = ad_0 = 0$$

Следовательно, $\mathrm{Im}\, D$ является $(\mathrm{Re}\, D)$-модулем.

Последовательность модулей

$$0 \longrightarrow \mathrm{Re}\, D \xrightarrow{\mathrm{id}} D \xrightarrow{\mathrm{Im}} \mathrm{Im}\, D \longrightarrow 0$$

является точной последовательностью. Согласно определению (8.4) отображения $\mathrm{Re}$, следующая диаграмма коммутативна

$$\begin{array}{ccc} \mathrm{Re}\, D & \xrightarrow{\mathrm{id}} & D \\ {\scriptstyle \mathrm{id}}\downarrow & \swarrow{\scriptstyle \mathrm{Re}} & \\ \mathrm{Re}\, D & & \end{array}$$

Согласно утверждению (2) предложения [1]-III.3.3,

$$D = \mathrm{id}(\mathrm{Re}\, D) \oplus \ker \mathrm{Re} \tag{8.7}$$

Согласно определению (8.5)

$$\ker \mathrm{Re} = \{d \in D : \mathrm{Re}\, d = 0\} = \mathrm{Im}\, D \tag{8.8}$$

Равенство (8.6) следует из равенств (8.7), (8.8). □

Согласно теореме 8.3, однозначно определенно представление

$$d = \mathrm{Re}\, d + \mathrm{Im}\, d \tag{8.9}$$

**Определение 8.4.** Отображение

$$d^* = \mathrm{Re}\, d - \mathrm{Im}\, d \tag{8.10}$$

называется **сопряжением в кольце** при условии, если это отображение удовлетворяет равенству

$$(cd)^* = d^*\, c^* \tag{8.11}$$

Кольцо $D$, в котором определено сопряжение, называется **кольцом с сопряжением**. □

Согласно теореме 8.1, кольцо с сопряжением $D$ является ассоциативной $(\mathrm{Re}\, D)$-алгеброй, которую мы также будем называть **алгеброй с сопряжением**.

**Следствие 8.5.**

$$(d^*)^* = d \tag{8.12}$$

□

Очевидно, существуют кольца, в которых условие (8.11) не выполнено. Для того, чтобы понять какие свойства у кольца с сопряжением $D$, рассмотрим структурные константы $(\mathrm{Re}\, D)$-алгебры $D$. Равенство

$$C^{l}_{0k} = C^{l}_{k0} = \delta^k_l \tag{8.13}$$

является следствием равенства $1d = d1 = d$.

**Теорема 8.6.** *Кольцо $D$ является кольцом с сопряжением, если*

$$C^{0}_{kl}=C^{0}_{lk}\quad C^{p}_{kl}=-C^{p}_{lk} \tag{8.14}$$

$$1<k<n\quad 1<l<n\quad 1<p<n$$

*Доказательство.* Из равенств (8.9), (8.10) следует

$$\begin{aligned}(cd)^*&=(\operatorname{Re}c\operatorname{Re}d+\operatorname{Re}c\operatorname{Im}d+\operatorname{Im}c\operatorname{Re}d+\operatorname{Im}c\operatorname{Im}d)^*\\&=\operatorname{Re}c\operatorname{Re}d-\operatorname{Re}c\operatorname{Im}d-\operatorname{Im}c\operatorname{Re}d+(\operatorname{Im}c\operatorname{Im}d)^*\end{aligned} \tag{8.15}$$

$$\begin{aligned}d^*\,c^*&=(\operatorname{Re}d-\operatorname{Im}d)(\operatorname{Re}c-\operatorname{Im}c)\\&=\operatorname{Re}d\operatorname{Re}c-\operatorname{Re}d\operatorname{Im}c-\operatorname{Im}d\operatorname{Re}c+\operatorname{Im}d\operatorname{Im}c\end{aligned} \tag{8.16}$$

Следовательно, из равенства (8.11) следует

$$(\operatorname{Im}c\operatorname{Im}d)^*=\operatorname{Im}d\operatorname{Im}c \tag{8.17}$$

Пусть $\overline{\overline{e}}$ - базис $(\operatorname{Re}D)$-алгебры $D$. Из равенства (8.17) следует

$$\begin{aligned}(C^{0}_{kl}c^kd^l+C^{p}_{kl}c^kd^l\overline{e}_p)^*&=C^{0}_{kl}c^kd^l-C^{p}_{kl}c^kd^l\overline{e}_p\\&=C^{0}_{kl}d^kc^l+C^{p}_{kl}d^kc^l\overline{e}_p\end{aligned} \tag{8.18}$$

$$1<k<n\quad 1<l<n\quad 1<p<n$$

Равенство (8.14) следует из равенства (8.18). □

**Следствие 8.7.**

$$\overline{e}_k\,\overline{e}_k\in\operatorname{Re}D$$

□

Отображение сопряжения можно представить с помощью матрицы $I$

$$\begin{aligned}d^*=I\circ d\qquad I^{0}_{k}&=\ \delta^{0}_{k}\quad\ I^{k}_{0}=\ \delta^{k}_{0}\qquad k=0,...,n\\I^{m}_{k}&=-\delta^{m}_{k}\quad I^{k}_{m}=-\delta^{k}_{m}\quad m=1,...,n\end{aligned} \tag{8.19}$$

**Пример 8.8.** Произведение в множестве комплексных чисел $C$ коммутативно. Однако поле комлексных чисел содержит подполе $R$. Векторное пространство $\operatorname{Im}C$ имеет размерность 1 и базис $\overline{e}_1=i$. Соответственно

$$c^*=c^0-c^1i$$

□

**Пример 8.9.** Тело кватернионов $H$ содержит подполе $R$. Векторное пространство $\operatorname{Im}H$ имеет размерность 3 и базис

$$\overline{e}_1=i\quad\overline{e}_2=j\quad\overline{e}_3=k$$

Соответственно

$$d^*=d^0-d^1i-d^2j-d^3k$$

□

## 9. Линейное отображение алгебры с сопряжением

Мы можем рассматривать кольцо $D$ как $D\star$-модуль размерности 1.

**Теорема 9.1.** *$D\star$-линейное отображение кольца с сопряжением $D$ имеет вид*

$$f(c)=cd\quad d\in D \tag{9.1}$$

*Пусть $\overline{\overline{e}}$ - $(\mathrm{Re}\,D)$-базис кольца $D$. Координаты отображения $f$*

$$f(c)=c^l f_l^i \overline{e}_i \tag{9.2}$$

*относительно базиса $\overline{\overline{e}}$ удовлетворяют равенству*

$$f_k^j=f_0^i C_{ki}^j \tag{9.3}$$

$$d=f_0^i \overline{e}_i \tag{9.4}$$

*Доказательство.* $(\mathrm{Re}\,D)$-линейное отображение удовлетворяет равенствам

$$f(ac)=(ac)^i f_i^j \overline{e}_j=a^k c^l C_{kl}^i f_i^j \overline{e}_j \tag{9.5}$$

$$af(c)=a^k (f(c))^i C_{ki}^j \overline{e}_j=a^k c^l f_l^i C_{ki}^j \overline{e}_j \tag{9.6}$$

Из равенств (9.5), (9.6), следует

$$a^k c^l C_{kl}^i f_i^j \overline{e}_j=a^k c^l f_l^i C_{ki}^j \overline{e}_j \tag{9.7}$$

Равенство

$$C_{kl}^i f_i^j=f_l^i C_{ki}^j \tag{9.8}$$

следует из равенства (9.7).

Из равенства (9.8) следует[6]

$$C_{k0}^i f_i^j=f_0^i C_{ki}^j \tag{9.9}$$

Из равенств (8.13), (9.9), следует

$$\delta_k^i f_i^j=f_0^i C_{ki}^j \tag{9.10}$$

Равенство (9.3) следует из равенства (9.10).

Из равенств (9.2), (9.3) следует

$$f(c)=c^l f_0^j C_{lj}^i \overline{e}_i \tag{9.11}$$

Равенство (9.4) следует из равенств (9.1), (9.11). □

**Определение 9.2.** Пусть $D$ - кольцо с сопряжением. Отображение

$$f:D\to D$$

называется **$D\star$-антилинейным**[7] если отображение $f$ удовлетворяет равенству

$$f(da)=f(a)d^* \tag{9.12}$$

□

[6]Пусть $l=0$.

[7]Определение 9.2 дано по аналогии с определением 5.1. Это определение также учитывает требования к антилинейному гомоморфизму в определении 10.2.

**Теорема 9.3.** *$D\star$-антилинейное отображение кольца с сопряжением $D$ имеет вид*

$$f(c)=dc^*=dI\circ c\quad d\in D \tag{9.13}$$

*Пусть $\overline{\overline{e}}$ - $(\mathrm{Re}\,D)$-базис кольца $D$. Координаты отображения $f$*

$$f(c)=c^l f^i_l \overline{e}_i \tag{9.14}$$

*относительно базиса $\overline{\overline{e}}$ удовлетворяют равенству*

$$f^j_l=f^k_0 I^i_l C^j_{ki} \tag{9.15}$$

$$d=f^i_0\overline{e}_i \tag{9.16}$$

*Доказательство.* $(\mathrm{Re}\,D)$-линейное отображение удовлетворяет равенствам

$$f(ac)=((ac)^*)^i f^j_i\overline{e}_j=c^p a^q (C^i_{pq})^* f^j_i\overline{e}_j \tag{9.17}$$

$$\begin{aligned}f(c)a^*&=(f(c))^k(a^*)^i C^j_{ki}\overline{e}_j=(c^*)^l f^k_l(a^*)^i C^j_{ki}\overline{e}_j\\&=I^l_p c^p f^k_l I^i_q a^q C^j_{ki}\overline{e}_j\end{aligned} \tag{9.18}$$

Из равенств (9.17), (9.18), следует

$$c^p a^q (C^i_{pq})^* f^j_i\overline{e}_j=I^l_p c^p f^k_l I^i_q a^q C^j_{ki}\overline{e}_j \tag{9.19}$$

Равенство

$$(C^i_{pq})^* f^j_i=I^l_p f^k_l I^i_q C^j_{ki} \tag{9.20}$$

следует из равенства (9.19).

Из равенства (9.20) следует[8]

$$(C^i_{0q})^* f^j_i=I^l_0 f^k_l I^i_q C^j_{ki} \tag{9.21}$$

Из равенств (8.13), (8.19), (9.21), следует

$$\delta^i_q f^j_i=\delta^l_0 f^k_l I^i_q C^j_{ki} \tag{9.22}$$

Равенство (9.15) следует из равенства (9.22).

Из равенств (9.14), (9.15) следует

$$f(c)=c^l f^k_0 I^i_l C^j_{ki}\overline{e}_j=f^k_0 I^i_l c^l C^j_{ki}\overline{e}_j=f^k_0 (c^*)^i C^j_{ki}\overline{e}_j \tag{9.23}$$

Равенство (9.16) следует из равенств (9.1), (9.23). □

**Теорема 9.4.** *Если отображение*

$$f:D\to D$$

*является $D\star$-линейным отображением, то отображение*

$$g:D\to D\quad g(d)=(f(d))^*$$

*является $D\star$-антилинейным отображением.*

*Доказательство.* Утверждение теоремы следует из равенства

$$g(cd)=(f(cd))^*=(cf(d))^*=(f(d))^*\,c^*=g(d)\,c^*$$

и определения 9.2. □

[8] Пусть $p=0$.

## 10. Антилинейный гомоморфизм $D\star$-алгебры

В этом разделе $D$ - это кольцо с сопряжением.

**Теорема 10.1.** *Пусть*

$$\begin{array}{ccccc} D & \xrightarrow[*]{f_{1\cdot 1,2}} & A_1 & \xrightarrow[*]{f_{1\cdot 2,3}} & A_1 \\ &&&& \uparrow {\scriptstyle *}\, f_{1\cdot 1,2} \\ &&&& D \end{array} \qquad \begin{array}{r} f_{1\cdot 1,2}(d):\overline{v}\to d\,\overline{v} \\ f_{1\cdot 2,3}(\overline{v}):\overline{w}\to \overline{C}_1(\overline{v},\overline{w}) \\ \overline{C}_1\in\mathcal{L}(A_1^2;A_1) \end{array} \tag{10.1}$$

*диаграмма представлений, описывающая $D\star$-алгебру $A_1$. Пусть*

$$\begin{array}{ccccc} D & \xrightarrow[*]{f_{2\cdot 1,2}} & A_2 & \xrightarrow[*]{f_{2\cdot 2,3}} & A_2 \\ &&&& \uparrow {\scriptstyle *}\, f_{2\cdot 1,2} \\ &&&& D \end{array} \qquad \begin{array}{r} f_{2\cdot 1,2}(d):\overline{v}\to \overline{v}\,d \\ f_{2\cdot 2,3}(\overline{v}):\overline{w}\to \overline{C}_2(\overline{w},\overline{v}) \\ \overline{C}_2\in\mathcal{L}(A_2^2;A_2) \end{array} \tag{10.2}$$

*диаграмма представлений, описывающая $\star D$-алгебру $A_2$. Морфизм $D\star$-алгебры $A_1$ в $\star D$-алгебру $A_2$*

$$r_1: D\to D \quad \overline{r}_2: A_1\to A_2$$
$$r_1(d)=d^*$$

*является $D\star$-антилинейным отображением $D\star$-алгебры $A_1$ в $\star D$-алгебру $A_2$ таким, что*

$$\overline{r}_2(\overline{a}\overline{b})=\overline{r}_2(\overline{b})\overline{r}_2(\overline{a}) \tag{10.3}$$

*Доказательство.* Согласно равенству [5]-(4.2.3), морфизм $(r_1,\overline{r}_2)$ представления $f_{1,2}$ удовлетворяет равенству

$$\begin{aligned} \overline{r}_2(f_{1\cdot 1,2}(d)(\overline{a})) &= f_{2\cdot 1,2}(r_1(d))(\overline{r}_2(\overline{a})) \\ \overline{r}_2(d\,\overline{a}) &= \overline{r}_2(\overline{a})\,d^* \end{aligned} \tag{10.4}$$

Следовательно, отображение $(r_1,\overline{r}_2)$ является $D\star$-антилинейным отображением.

Согласно равенству [5]-(4.2.3), морфизм $(\overline{r}_2,\overline{r}_2)$ представления $f_{2,3}$ удовлетворяет равенству[9]

$$\overline{r}_2(f_{1\cdot 2,3}(\overline{a}_2)(\overline{a}_3)) = f_{2\cdot 2,3}(\overline{r}_2(\overline{a}_2))(\overline{r}_2(\overline{a}_3)) \tag{10.5}$$

Из равенств (10.5), (10.1), (10.2), следует

$$\overline{r}_2(\overline{C}_1(\overline{v},\overline{w}))=\overline{C}_2(\overline{r}_2(\overline{w}),\overline{r}_2(\overline{v})) \tag{10.6}$$

Равенство (10.3) следует из равенств (10.6), (6.3). □

**Определение 10.2.** Морфизм представлений $D\star$-алгебры $A_1$ в $\star D$-алгебру $A_2$ называется **$D\star$-антилинейным гомоморфизмом** $D\star$-алгебры $A_1$ в $\star D$-алгебру $A_2$. □

[9]Так как в диаграммах представлений (6.4), (6.5), носители $\Omega_2$-алгебры и $\Omega_3$-алгебры совпадают, то также совпадают морфизмы представлений на уровнях 2 и 3.

**Теорема 10.3.** *Пусть $\overline{\overline{e}}_1$ - базис $D\star$-алгебры $A_1$. Пусть $\overline{\overline{e}}_2$ - базис $\star D$-алгебры $A_2$. Тогда антилинейный гомоморфизм*[10] $(r_1, \overline{r}_2)$ *$D\star$-алгебры $A_1$ в $\star D$-алгебру $A_2$ имеет представление*[11]

$$\overline{b} = \overline{e}_{2*}{}^* r_{2*}{}^* r_1(a) = \overline{e}_{2\cdot i}\, r_{2\cdot}{}^{i}_{j}\, r_1(a^j) \tag{10.7}$$
$$b = r_{2*}{}^* r_1(a)$$

*относительно заданных базисов. Здесь*

- *$a$ - координатная матрица вектора $\overline{a}$ относительно базиса $\overline{\overline{e}}_1$.*
- *$b$ - координатная матрица вектора*
  $$\overline{b} = \overline{r}_2(\overline{a})$$
  *относительно базиса $\overline{\overline{e}}_2$.*
- *$r_2$ - координатная матрица множества векторов $(\overline{r}_2(\overline{e}_{1\cdot i}))$ относительно базиса $\overline{\overline{e}}_2$. Мы будем называть матрицу $r_2$* **матрицей антилинейного гомоморфизма** *относительно базисов $\overline{\overline{e}}_1$ и $\overline{\overline{e}}_2$.*

*Доказательство.* Вектор $\overline{a} \in A_1$ имеет разложение
$$\overline{a} = a^*{}_* \overline{e}_1$$
относительно базиса $\overline{\overline{e}}_1$. Вектор $\overline{b} \in A_2$ имеет разложение
$$\overline{b} = \overline{e}_{2*}{}^* b \tag{10.8}$$
относительно базиса $\overline{\overline{e}}_2$.

Так как $(r_1, \overline{r}_2)$ - антилинейный гомоморфизм, то на основании (10.4) следует, что
$$\overline{b} = \overline{r}_2(\overline{a}) = \overline{r}_2(a^*{}_* \overline{e}_1) = \overline{r}_2(\overline{e}_1)_*{}^* a^* \tag{10.9}$$
где
$$a^* = \begin{pmatrix} (a^1)^* \\ ... \\ (a^n)^* \end{pmatrix}$$
$\overline{r}_2(\overline{e}_{1\cdot i})$ также вектор $\star D$-модуля $A_2$ и имеет разложение
$$\overline{r}_2(\overline{e}_{1\cdot i}) = \overline{e}_{2*}{}^* r_{2\cdot i} = \overline{e}_j\, r_{2\cdot}{}^{j}_{i} \tag{10.10}$$
относительно базиса $\overline{\overline{e}}_2$. Комбинируя (10.9) и (10.10), мы получаем
$$\overline{b} = \overline{e}_{2\circ}{}^\circ r_{2\circ}{}^\circ a^* \tag{10.11}$$
(10.7) следует из сравнения (10.8) и (10.11) и теоремы [2]-4.3.3. □

**Теорема 10.4.** *Пусть $\overline{\overline{e}}_1$ - базис $D\star$-алгебры $A_1$. Пусть $\overline{\overline{e}}_2$ - базис $\star D$-алгебры $A_2$. Матрица антилинейного гомоморфизма и структурные константы связаны соотношением*
$$(C_{1\cdot}{}^{k}_{ij})^*\, r_{2\cdot}{}^{l}_{k} = r_{2\cdot}{}^{p}_{j}\, r_{2\cdot}{}^{q}_{i}\, C_{2\cdot}{}^{l}_{pq} \tag{10.12}$$

---

[10] Эта теорема аналогична теореме [2]-4.4.3.

[11] Для поля комплексных чисел, равенство (5.11) является следствием равенства (10.7).

*Доказательство.* Пусть

$$\overline{a},\overline{b}\in A_1\quad \overline{a}=a^*{}_*\overline{e}_1\quad \overline{b}=b^*{}_*\overline{e}_1$$

Из равенств (6.2), (6.3), (10.1), следует

$$\overline{a}\overline{b}=a^ib^jC_{1\cdot}{}^k_{ij}\overline{e}_{1\cdot k}\tag{10.13}$$

Из равенств (10.4), (10.13), следует

$$\overline{r}_2(\overline{a}\overline{b})=\overline{r}_2(\overline{e}_{1\cdot k})(a^ib^jC_{1\cdot}{}^k_{ij})^*=\overline{r}_2(\overline{e}_{1\cdot k})(C_{1\cdot}{}^k_{ij})^*(b^j)^*(a^i)^*\tag{10.14}$$

Из теоремы 10.3 и равенства (10.14), следует

$$\overline{r}_2(\overline{a}\overline{b})=\overline{e}_{2\cdot l}\,r_{2\cdot}{}^l_k(C_{1\cdot}{}^k_{ij})^*(a^i)^*(b^j)^*\tag{10.15}$$

Из равенства (10.3) и теоремы 10.3, следует

$$\overline{r}_2(\overline{a}\overline{b})=\overline{r}_2(\overline{b})\overline{r}_2(\overline{a})=\overline{e}_{2\cdot q}\,r_{2\cdot}{}^q_j(b^j)^*\,\overline{e}_{2\cdot p}\,r_{2\cdot}{}^p_i(a^i)^*\tag{10.16}$$

Из равенств (6.2), (6.3), (10.2), (10.16), следует

$$\overline{r}_2(\overline{a}\overline{b})=\overline{e}_{2\cdot l}\,C_{2\cdot}{}^l_{pq}r_{2\cdot}{}^p_i(a^i)^*r_{2\cdot}{}^q_j(b^j)^*\tag{10.17}$$

Из равенств (10.15), (10.17), следует

$$\overline{e}_{2\cdot l}\,r_{2\cdot}{}^l_k(C_{1\cdot}{}^k_{ij})^*(a^i)^*(b^j)^*=\overline{e}_{2\cdot l}\,C_{2\cdot}{}^l_{pq}r_{2\cdot}{}^p_i(a^i)^*r_{2\cdot}{}^q_j(b^j)^*\tag{10.18}$$

Равенство (10.12) следует из равенства (10.18), так как векторы базиса $\overline{\overline{e}}_2$ линейно независимы, и $a^i$, $b^i$ (а следовательно, $(a^i)^*$, $(b^i)^*$) - произвольные величины. □

## 11. Инволюция банаховой $D$-алгебры

В этом разделе $D$ - это кольцо с сопряжением. Мы будем также полагать, что в абелевой группе $A$ определены структуры $D\star$-модуля и $\star D$-модуля.[12]

**Определение 11.1.** Пусть $A$ - банаховая $D$-алгебра. Инволюция - это антилинейный гомоморфизм $x\to x^*$, $x\in A$, такое, что

$$x^{**}=x\tag{11.1}$$

$$\|x^*x\|=\|x\|^2$$

$$x,y\in A$$

□

**Теорема 11.2.** *Пусть $A$ - банаховая $D$-алгебра. Пусть отображение $x\to x^*$, $x\in A$, является инволюцией. Тогда*

$$(xy)^*=y^*x^*\tag{11.2}$$

*Доказательство.* Равенство (11.2) является следствием определения 11.1 и теоремы 10.1. □

**Определение 11.3.** Базис $\overline{\overline{e}}$ называется нормированным базисом, если $\|\overline{e}_i\|=1$ для любого вектора $\overline{e}_i$ базиса $\overline{\overline{e}}$. □

[12]Условия, когда данное требование выполнено, будут рассмотренны в отдельной статье. Если $D$ - поле, то это условие очевидно. Если $D$ - тело, то это условие выполнено согласно теореме [2]-6.3.1.

Не нарушая общности, мы можем положить, что базис $\overline{\overline{e}}$ нормирован. Если предположить, что норма вектора $\overline{e}_i$ отлична от 1, то мы можем этот вектор заменить вектором

$$\overline{e}'_i = \frac{1}{\|\overline{e}_i\|}\overline{e}_i$$

**Теорема 11.4.** *Пусть $D$ - коммутативное кольцо. Если отображение*

$$f : A_1 \to A_2$$

*является линейным, то отображение*

$$h(x) = f(x^*)$$

*является антилинейным. Если отображение*

$$f : A_1 \to A_2$$

*является линейным гомоморфизмом, то отображение*

$$h(x) = f(x^*)$$

*является антилинейным гомоморфизмом.*

*Доказательство.* Согласно теореме 10.1 и равенству (10.4)

$$(da)^* = a^*\,d^* \quad d \in D \quad a \in A \tag{11.3}$$

Если $f$ - линейное отображение, то

$$f(dx) = df(x) \tag{11.4}$$

Из равенств (11.3), (11.4), следует

$$h(xd) = f((xd)^*) = f(d^*\,x^*) = d^*\,f(x^*) = d^*\,h(x)$$

Если $f$ - линейный гомоморфизм, то из равенств (6.6), (11.2), следует

$$h(ab) = f((ab)^*) = f(b^*\,a^*) = f(b^*)f(a^*) = h(b)h(a)$$

Следовательно, $h$ является антилинейным гомоморфизмом. □

## 12. $C^*$-алгебра

**Определение 12.1.** **$C^*$-алгебра** $A$ - это банаховая $C$-алгебра, в которой определена инволюция.[13] □

Согласно теореме 5.4, инволюция имеет вид

$$\begin{aligned}(x^*)^{0j} &= f_{0i}^{0j}x^{0i} + f_{0i}^{1j}x^{1i}\\(x^*)^{1j} &= f_{0i}^{1j}x^{0i} - f_{0i}^{0j}x^{1i}\end{aligned} \tag{12.1}$$

относительно базиса $\overline{\overline{e}}_{AR}$ .

**Теорема 12.2.**

$$\begin{aligned}\delta_k^j &= f_{0i}^{0j}f_{0k}^{0i} + f_{0i}^{1j}f_{0k}^{1i}\\0 &= f_{0i}^{0j}f_{0k}^{1i} - f_{0i}^{1j}f_{0k}^{0i}\end{aligned} \tag{12.2}$$

---

[13] Это определение сделано на основе определения [6]-2.2.1.

*Доказательство.* Из равенства (12.1) следует

$$
\begin{aligned}
(x^{**})^{0j} &= f^{0j}_{0i}(x^*)^{0i} + f^{1j}_{0i}(x^*)^{1i}\\
&= f^{0j}_{0i}(f^{0i}_{0k}x^{0k} + f^{1i}_{0k}x^{1k}) + f^{1j}_{0i}(f^{1i}_{0k}x^{0k} - f^{0i}_{0k}x^{1k})\\
&= f^{0j}_{0i}f^{0i}_{0k}x^{0k} + f^{0j}_{0i}f^{1i}_{0k}x^{1k} + f^{1j}_{0i}f^{1i}_{0k}x^{0k} - f^{1j}_{0i}f^{0i}_{0k}x^{1k}\\
(x^{**})^{1j} &= f^{1j}_{0i}(x^*)^{0i} - f^{0j}_{0i}(x^*)^{1i}\\
&= f^{1j}_{0i}(f^{0i}_{0k}x^{0k} + f^{1i}_{0k}x^{1k}) - f^{0j}_{0i}(f^{1i}_{0k}x^{0k} - f^{0i}_{0k}x^{1k})\\
&= f^{1j}_{0i}f^{0i}_{0k}x^{0k} + f^{1j}_{0i}f^{1i}_{0k}x^{1k} - f^{0j}_{0i}f^{1i}_{0k}x^{0k} + f^{0j}_{0i}f^{0i}_{0k}x^{1k}
\end{aligned}
\tag{12.3}
$$

(12.2) следует из равенств (11.1), (12.3). □

**Теорема 12.3.**

$$
\begin{aligned}
&C_{AC\cdot}{}^{k0}_{ij}f^{0l}_{0k} + C_{AC\cdot}{}^{k1}_{ij}f^{1l}_{0k}\\
&= f^{0p}_{0j}f^{0q}_{0i}C_{AC\cdot}{}^{l0}_{pq} - f^{1p}_{0j}f^{1q}_{0i}C_{AC\cdot}{}^{l0}_{pq} - f^{0p}_{0j}f^{1q}_{0i}C_{AC\cdot}{}^{l1}_{pq} - f^{1p}_{0j}f^{0q}_{0i}C_{AC\cdot}{}^{l1}_{pq}\\
&C_{AC\cdot}{}^{k0}_{ij}f^{1l}_{0k} - C_{AC\cdot}{}^{k1}_{ij}f^{0l}_{0k}\\
&= f^{0p}_{0j}f^{0q}_{0i}C_{AC\cdot}{}^{l1}_{pq} - f^{1p}_{0j}f^{1q}_{0i}C_{AC\cdot}{}^{l1}_{pq} + f^{0p}_{0j}f^{1q}_{0i}C_{AC\cdot}{}^{l0}_{pq} + f^{1p}_{0j}f^{0q}_{0i}C_{AC\cdot}{}^{l0}_{pq}
\end{aligned}
\tag{12.4}
$$

*Доказательство.* Из равенства (5.11) следует

$$
f^j_i = f^{0j}_{0i} + f^{1j}_{0i}i
\tag{12.5}
$$

Из равенств (3.17), (12.5), (10.12), следует

$$
\begin{aligned}
&(C_{AC\cdot}{}^k_{ij})^*f^l_k = f^p_j f^q_i C_{AC\cdot}{}^l_{pq}\\
&(C_{AC\cdot}{}^{k0}_{ij} + C_{AC\cdot}{}^{k1}_{ij}i)^*(f^{0l}_{0k} + f^{1l}_{0k}i)\\
&= (f^{0p}_{0j} + f^{1p}_{0j}i)(f^{0q}_{0i} + f^{1q}_{0i}i)(C_{AC\cdot}{}^{l0}_{pq} + C_{AC\cdot}{}^{l1}_{pq}i)
\end{aligned}
\tag{12.6}
$$

Из равенства (12.6) следует

$$
\begin{aligned}
&C_{AC\cdot}{}^{k0}_{ij}f^{0l}_{0k} + C_{AC\cdot}{}^{k1}_{ij}f^{1l}_{0k} + (C_{AC\cdot}{}^{k0}_{ij}f^{1l}_{0k} - C_{AC\cdot}{}^{k1}_{ij}f^{0l}_{0k})i\\
&= (f^{0p}_{0j}f^{0q}_{0i} - f^{1p}_{0j}f^{1q}_{0i} + (f^{0p}_{0j}f^{1q}_{0i} + f^{1p}_{0j}f^{0q}_{0i})i)(C_{AC\cdot}{}^{l0}_{pq} + C_{AC\cdot}{}^{l1}_{pq}i)\\
&= f^{0p}_{0j}f^{0q}_{0i}C_{AC\cdot}{}^{l0}_{pq} - f^{1p}_{0j}f^{1q}_{0i}C_{AC\cdot}{}^{l0}_{pq} - f^{0p}_{0j}f^{1q}_{0i}C_{AC\cdot}{}^{l1}_{pq} - f^{1p}_{0j}f^{0q}_{0i}C_{AC\cdot}{}^{l1}_{pq}\\
&+ (f^{0p}_{0j}f^{0q}_{0i}C_{AC\cdot}{}^{l1}_{pq} - f^{1p}_{0j}f^{1q}_{0i}C_{AC\cdot}{}^{l1}_{pq} + f^{0p}_{0j}f^{1q}_{0i}C_{AC\cdot}{}^{l0}_{pq} + f^{1p}_{0j}f^{0q}_{0i}C_{AC\cdot}{}^{l0}_{pq})i
\end{aligned}
\tag{12.7}
$$

Равенство (12.4) следует из равенства (12.7). □

## 13. Пример инволюции

Пусть

$$
e^1_1 = \begin{pmatrix} 1 & 0\\ 0 & 0 \end{pmatrix} \quad e^2_1 = \begin{pmatrix} 0 & 1\\ 0 & 0 \end{pmatrix}
$$

$$
e^1_2 = \begin{pmatrix} 0 & 0\\ 1 & 0 \end{pmatrix} \quad e^2_2 = \begin{pmatrix} 0 & 0\\ 0 & 1 \end{pmatrix}
$$

базис $\overline{\overline{e}}$ алгебры $2 \times 2$ матриц. Матрица

$$a = \begin{pmatrix} a_1^1 & a_2^1 \\ a_1^2 & a_2^2 \end{pmatrix}$$

имеет разложение

$$a = a_1^1 e_1^1 + a_2^1 e_1^2 + a_1^2 e_2^1 + a_2^2 e_2^2$$

относительно базиса $\overline{\overline{e}}$. Структурные константы в алгебре матриц имеют вид

$$C^{\cdot p}_{q\cdot}{}^{i}_{j\cdot}{}^{k}_{l} = \delta^i_l \delta^p_j \delta^k_q$$

Произведение матриц $a$ и $b$ имеет вид

$$(ab)^p_q = C^{\cdot p}_{q\cdot}{}^{i}_{j\cdot}{}^{k}_{l} a^j_i b^l_k = \delta^i_l \delta^p_j \delta^k_q a^j_i b^l_k = a^p_l b^l_q$$

Если мы матрицу $a$ представим в виде вектора

$$a = \begin{pmatrix} a_1^1 \\ a_2^1 \\ a_1^2 \\ a_2^2 \end{pmatrix}$$

то матрицу преобразование $f$

$$a'^i_j = f^{\cdot i}_{j\cdot}{}^{k}_{l} a^l_k$$

можно представить в следующем виде

$$f = \begin{pmatrix} f^{\cdot 1}_{1\cdot}{}^{1}_{1} & f^{\cdot 1}_{1\cdot}{}^{1}_{2} & f^{\cdot 1}_{1\cdot}{}^{2}_{1} & f^{\cdot 1}_{1\cdot}{}^{2}_{2} \\ f^{\cdot 1}_{2\cdot}{}^{1}_{1} & f^{\cdot 1}_{2\cdot}{}^{1}_{2} & f^{\cdot 1}_{2\cdot}{}^{2}_{1} & f^{\cdot 1}_{2\cdot}{}^{2}_{2} \\ f^{\cdot 2}_{1\cdot}{}^{1}_{1} & f^{\cdot 2}_{1\cdot}{}^{1}_{2} & f^{\cdot 2}_{1\cdot}{}^{2}_{1} & f^{\cdot 2}_{1\cdot}{}^{2}_{2} \\ f^{\cdot 2}_{2\cdot}{}^{1}_{1} & f^{\cdot 2}_{2\cdot}{}^{1}_{2} & f^{\cdot 2}_{2\cdot}{}^{2}_{1} & f^{\cdot 2}_{2\cdot}{}^{2}_{2} \end{pmatrix}$$

Преобразованию транспонирования

$$a = \begin{pmatrix} a_1^1 & a_2^1 \\ a_1^2 & a_2^2 \end{pmatrix} \to a^T = \begin{pmatrix} a_1^1 & a_1^2 \\ a_2^1 & a_2^2 \end{pmatrix}$$

соответствует матрица

$$f = \begin{pmatrix} 1 & 0 & 0 & 0 \\ 0 & 0 & 1 & 0 \\ 0 & 1 & 0 & 0 \\ 0 & 0 & 0 & 1 \end{pmatrix}$$

Так как мы не имеем аналитической записи для оператора транспонирования, то мы воспользуемся нестандартным форматом записи.

$$\begin{aligned}
(a^Tb^T)^p_q &= C^{\cdot p\ i\ k}_{q\cdot j\cdot l}(a^T)^j_i(b^T)^l_k\\
&= \delta^i_l\delta^p_j\delta^k_q\begin{pmatrix}1&0&0&0\\0&0&1&0\\0&1&0&0\\0&0&0&1\end{pmatrix}\begin{pmatrix}a^1_1\\a^1_2\\a^2_1\\a^2_2\end{pmatrix}\begin{pmatrix}1&0&0&0\\0&0&1&0\\0&1&0&0\\0&0&0&1\end{pmatrix}\begin{pmatrix}b^1_1\\b^1_2\\b^2_1\\b^2_2\end{pmatrix}\\
&= \delta^i_l\delta^p_j\delta^k_q\begin{pmatrix}{}^j_i:{}^1_1&a^1_1\\{}^j_i:{}^1_2&a^2_1\\{}^j_i:{}^2_1&a^1_2\\{}^j_i:{}^2_2&a^2_2\end{pmatrix}\begin{pmatrix}{}^l_k:{}^1_1&b^1_1\\{}^l_k:{}^1_2&b^2_1\\{}^l_k:{}^2_1&b^1_2\\{}^l_k:{}^2_2&b^2_2\end{pmatrix}=\begin{pmatrix}{}^p_q:{}^1_1&a^1_1b^1_1+a^2_1b^1_2\\{}^p_q:{}^1_2&a^1_1b^2_1+a^2_1b^2_2\\{}^p_q:{}^2_1&a^1_2b^1_1+a^2_2b^1_2\\{}^p_q:{}^2_2&a^1_2b^2_1+a^2_2b^2_2\end{pmatrix} \qquad (13.1)\\
&= \begin{pmatrix}1&0&0&0\\0&0&1&0\\0&1&0&0\\0&0&0&1\end{pmatrix}\begin{pmatrix}b^1_ka^k_1\\b^1_ka^k_2\\b^2_ka^k_1\\b^2_ka^k_2\end{pmatrix}\\
&= ((ba)^T)^p_q
\end{aligned}$$

Если $a^i_j$ - комплексные числа, то инволюцию можно представить антилинейным отображением с матрицей

$$f=\begin{pmatrix}I&0&0&0\\0&0&I&0\\0&I&0&0\\0&0&0&I\end{pmatrix}=\left(\begin{array}{cc:cc:cc:cc}1&0&0&0&0&0&0&0\\0&-1&0&0&0&0&0&0\\\hdashline 0&0&0&0&1&0&0&0\\0&0&0&0&0&-1&0&0\\\hdashline 0&0&1&0&0&0&0&0\\0&0&0&-1&0&0&0&0\\\hdashline 0&0&0&0&0&0&1&0\\0&0&0&0&0&0&0&-1\end{array}\right)$$

## 14. Список литературы

[1] Серж Ленг, Алгебра, М. Мир, 1968

[2] Александр Клейн, Лекции по линейной алгебре над телом, eprint arXiv:math.GM/0701238 (2010)

[3] Александр Клейн, Этюд о кватернионах, eprint arXiv:0909.0855 (2010)

[4] Aleks Kleyn, Linear Mappings of Free Algebra: First Steps in Noncommutative Linear Algebra, Lambert Academic Publishing, 2010

[5] Aleks Kleyn, Representation Theory: Representation of Universal Algebra, Lambert Academic Publishing, 2011
[6] William Arveson, A short course on spectral theory. Springer - Verlag, New York, 2002
[7] Robert Hermann, Topics in the mathematics of quantum mechanics. Math Sci Press, 1973

# 15. Предметный указатель

## 16. Специальные символы и обозначения

$d^*$ сопряжение в кольце 17

$\operatorname{Im} D$ модуль векторов кольца $D$ 16
$\operatorname{Im} d$ вектор элемента $d$ кольца 16

$\operatorname{Re} D$ алгебра скаляров кольца $D$ 16
$\operatorname{Re} d$ скаляр элемента $d$ кольца 16